%% file: Heat_Imperfect.tex
\documentclass[11pt]{article}

\usepackage{fullpage}
\usepackage{latexsym}
\usepackage{amsfonts}
\usepackage{amssymb}
\usepackage{amsmath}
\usepackage{subcaption}
\usepackage{color}
\usepackage{graphicx, graphics}
\usepackage{hyperref}
\usepackage{enumerate}
\usepackage{pstricks, pst-plot, epsfig}
\usepackage{xargs}
\usepackage{arydshln}
\newcommand{\ud}{\,\mathrm{d}}
\renewcommand\Re{\operatorname{Re}}
\renewcommand\Im{\operatorname{Im}}

\def\XXint#1#2#3{{\setbox0=\hbox{$#1{#2#3}{\int}$}
\vcenter{\hbox{$#2#3$}}\kern-.5\wd0}}

\definecolor{green}{RGB}{0,128,0}

\newcommand{\CC}{\mathbb{C}}

\begin{document}

\title{Multilayer diffusion in a composite medium with imperfect contact}
\date{\today}
\author{Natalie E. Sheils\\Institute for Mathematics and Its Applications\\nesheils@umn.edu}
\maketitle


\begin{abstract}
The problem of heat conduction in one-dimensional piecewise homogeneous composite materials is examined by providing an explicit solution of the one-dimensional heat equation in each domain. The location of the interfaces is known, but neither temperature nor heat flux are prescribed there.  We find a solution using the Unified Transform Method, due to Fokas and collaborators, applied to interface problems and compute solutions numerically.
\end{abstract}

\section{Introduction}
The problem of heat conduction in a composite wall is a classical problem in design and construction.  It is usual to restrict to the case of walls with physical properties that are constant throughout the material and are considered to be of infinite extent in the directions parallel to the wall. Further, we assume that temperature and heat flux do not vary in these directions. In that case, the mathematical model for heat conduction in each wall layer is given by~\cite[Chapter~10]{HahnO}:

\begin{subequations}\label{heatgeneral}
\begin{align}
u_t^{(j)}=\kappa_j u_{xx}^{(j)},  &&&x_{j-1}<x<x_j,\\
u^{(j)}(x,0)=u_0^{(j)}(x),  &&&x_{j-1}<x<x_j,\\
\beta_1 u^{(1)}(x_0,t)+\beta_2 u^{(1)}_x(x_0,t)=f_1(t),&&&t>0,\label{boundary1}\\
\beta_3 u^{(n+1)}(x_{n+1},t)+\beta_4 u^{(n+1)}_x(x_{n+1},t)=f_2(t), &&& t>0, \label{boundary2}
\end{align}
\end{subequations}
where $u^{(j)}(x,t)$ denotes the temperature in the wall layer indexed by $(j)$, $\kappa_j>0$ is the heat-conduction coefficient of the $j$-th layer (the inverse of its thermal diffusivity), $x=x_{j-1}$ is the left extent of the layer, $x=x_j$ is its right extent, and $\beta_n$ for $n=1,2,3,4$ are constants. The sub-indices denote derivatives with respect to the one-dimensional spatial variable $x$ and the temporal variable $t$. The function $u^{(j)}_0(x)$ is the prescribed initial condition of the system. The continuity of the temperature $u^{(j)}$ and of its associated heat flux $\kappa_j u^{(j)}_x$ are imposed across the interface between layers. In what follows it is convenient to use the quantity $\sigma_j$, defined as the positive square root of $\kappa_j$: $\sigma_j=\sqrt{\kappa_j}$.

If each layer is in perfect thermal contact then the interface conditions are
\begin{subequations}
\begin{align}
u^{(j)}(x_j,t)=u^{(j+1)}(x_j,t),  &&&t>0,\label{p_jump1}\\
\sigma_j^2 u^{(j)}_x(x_j,t)=\sigma_{j+1}^2u^{(j+1)}_x(x_j,t),  &&&t>0. \label{p_jump2}
\end{align}
\end{subequations}
A derivation of the interface conditions for perfect thermal contact is found in~\cite[Chapter~1]{HahnO}.  However, if the thermal contact is imperfect we prescribe the interface conditions
\begin{subequations}
\begin{align}
\sigma_j^2 u_x^{(j)}(x_j,t)=H_j\left(u^{(j+1)}(x_{j},t)-u^{(j)}(x_j,t)\right),  &&&t>0,\label{ip_jump1}\\
\sigma_{j+1}^2 u_x^{(j+1)}(x_j,t)=H_j\left(u^{(j+1)}(x_{j},t)-u^{(j)}(x_j,t)\right),  &&&t>0, \label{ip_jump2}
\end{align}
\end{subequations}
where $H_j\neq0$ is the contact transfer coefficient at $x=x_j$ and $1\leq j\leq n$.  Perfect thermal contact, is recovered in the limit $H_j\to\infty$.  In applications, imperfect boundary conditions are used to model roughness and contact resistance~\cite{BarrySweatman, CarrTurner, CarslawJaeger, SridharYovanovich}. Carr and Turner~\cite{CarrTurner} approach this problem using a semi-analytical method based on the Laplace transform and an orthogonal eigenfunction expansion.  Their interest in the problem is to accurately solve a two-scale modeling problem for transport or fluid flow in porous media exhibiting small scale heterogeneities in material properties.  The authors note that for a large number of layers, multilayer diffusion is possibly the most simple example of such a problem. However, their numerical implementation for their analytical solution only works for up to ten layers~\cite{CarrTurner}.   They also propose a ``semi-analytical" model which works for a large number of layers.

In this paper, we use the Fokas Method (also called the Unified Transform Method)~\cite{DeconinckTrogdonVasan, FokasBook, FokasPelloni4} to provide explicit solution formulae for different heat transport interface problems of the types described above.  Even for a simple problem (two finite walls in perfect thermal contact), the classical approach using separation of variables~\cite{HahnO} can provide an explicit answer only implicitly.  Indeed, the solution obtained in~\cite{HahnO} depends on certain eigenvalues defined through a transcendental equation that can be solved only numerically. In contrast, the Fokas Method produces an explicit solution formula involving only known quantities. In~\cite{DeconinckPelloniSheils} the problem of heat conduction in perfect thermal contact was considered using the Fokas Method to provide explicit solution formulae for a number of examples for up to three domains.  In this paper we extend that method to include more general interface conditions and a generic number of interfaces.

Interface problems for partial differential equations (PDEs) are initial boundary value problems for which the solution of an equation in one domain prescribes boundary conditions for the equations in adjacent domains. In applications, interface conditions are often obtained from conservation laws~\cite{Kevorkian}. Few interface problems allow for an explicit closed-form solution using classical solution methods. Using the Fokas Method, such solutions may be constructed for both dissipative and dispersive linear interface problems as shown in~\cite{Asvestas, DeconinckPelloniSheils,Mantzavinos, SheilsDeconinck_PeriodicHeat, SheilsDeconinck_LS, SheilsDeconinck_LSp, SheilsSmith}. 

\section{The Fokas Method for the heat equation}
We follow the standard steps in the Fokas Method.  Assuming existence of a solution, we begin with the so-called ``local relations'':

\begin{equation}\label{localrelation}
\left(e^{-ikx+\omega_j(k) t}u^{(j)}\right)_t=\left(e^{-ikx+\omega_j(k) t}\sigma_j^2( u^{(j)}_x+iku^{(j)})\right)_x,
\end{equation}
where $\omega_j(k)=(\sigma_jk)^2.$  Without loss of generality we shift the problem so that $x_0=0$.

\begin{figure}[htbp]
\begin{center}
\def\svgwidth{.75\textwidth}
   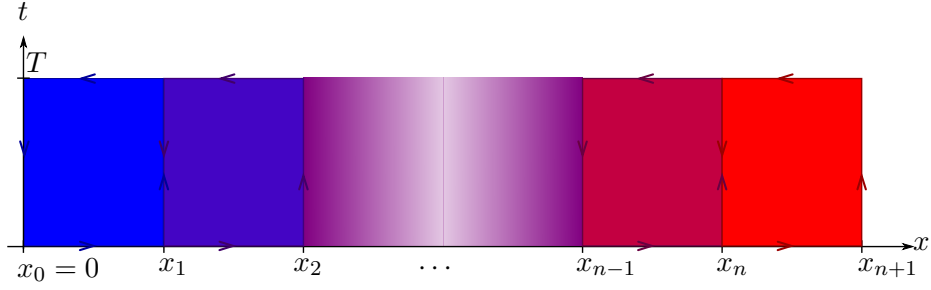 
   \caption{Domains for the application of Green's Theorem in the case of a finite domain with $n$ interfaces.   \label{fig:GR_domain_nf}}
  \end{center}
\end{figure}

Integrating each local relation~\eqref{localrelation} around the appropriate domain (see Figure~\ref{fig:GR_domain_nf}) and applying Green's Theorem we find the global relations:
\begin{equation}\label{nGR}
\begin{split}
0=&\int_{x_{j-1}}^{x_j} e^{-ikx}u^{(j)}_0(x)\ud x-\int_{x_{j-1}}^{x_j} e^{-ikx+\omega_j(k) T}u^{(j)}(x,T)\ud x\\
&+\int_0^T\sigma_j^2 e^{-ikx_j+\omega_j(k) s} (u^{(j)}_x(x_j,s)+ i k u^{(j)}(x_j,s))\ud s\\
&-\int_0^T\sigma_j^2 e^{-ikx_{j-1}+\omega_j(k) s} (u^{(j)}_x(x_{j-1},s)+ i k u^{(j)}(x_{j-1},s))\ud s,
\end{split}
\end{equation}
for $1\leq j\leq n+1$.  

We define the following transforms for $1\leq j\leq n+1:$
\begin{align*}
\hat{u}^{(j)}(k,t)&=\int_{x_{j-1}}^{x_j} e^{-ikx} u^{(j)}(x,t)\ud x,&x_{j-1}<x<x_j, &~~t>0,\\
\hat{u}^{(j)}_0(k)&=\int_{x_{j-1}}^{x_j} e^{-ikx} u_0^{(j)}(x)\ud x,&x_{j-1}<x<x_j,\\
g_0^{(j)}(\omega,t)&=\int_0^t e^{\omega s} u^{(j)}(x_{j-1},s)\ud s,&&~~t>0,\\
g_1^{(j)}(\omega,t)&=\int_0^t e^{\omega s} u^{(j)}_x(x_{j-1},s)\ud s,&&~~t>0,\\
h_0^{(j)}(\omega,t)&=\int_0^t e^{\omega s} u^{(j)}(x_j,s)\ud s,&&~~t>0,\\
h_1^{(j)}(\omega,t)&=\int_0^t e^{\omega s} u^{(j)}_x(x_j,s)\ud s,&&~~t>0.
\end{align*}
All of these integrals are proper integrals defined for $k\in\CC$.  With these definitions the global relations become
\begin{equation}\label{GR}
\begin{split}
e^{\omega_j(k) T}\hat{u}^{(j)}(k,T)=&\hat{u}^{(j)}_0(k)+\sigma_j^2 e^{-ikx_j} \left(h_1^{(j)}(\omega_j(k),T)+ i k h_0^{(j)}(\omega_j(k),T)\right)\\
&-\sigma_j^2 e^{-ikx_{j-1}} \left(g_1^{(j)}(\omega_j(k),T)+ i k g_0^{(j)}(\omega_j(k),T)\right),
\end{split}
\end{equation}
for $1\leq j\leq n+1$.  We transform the global relations so that $g_\ell^{(j)}(\cdot,T)$ and $h_\ell^{(j)}(\cdot,T)$ depend on a common argument  $\nu^2$ through the change of variables $k=\nu/\sigma_j$:

\begin{equation}\label{GR_common_plus}
\begin{split}
e^{\nu^2 T}\hat{u}^{(j)}\left(\frac{\nu}{\sigma_j},T\right)=&\hat{u}^{(j)}_0\left(\frac{\nu}{\sigma_j}\right)+e^{-\frac{i\nu x_j}{\sigma_j}} \left(\sigma_j^2 h_1^{(j)}(\nu^2,T)+ i \sigma_j\nu h_0^{(j)}(\nu^2,T)\right)\\
&-e^{-\frac{i\nu x_{j-1}}{\sigma_j}} \left(\sigma_j^2 g_1^{(j)}(\nu^2,T)+ i \sigma_j\nu g_0^{(j)}(\nu^2,T)\right),
\end{split}
\end{equation}
where $1\leq j\leq n+1$.  

The dispersion relation $\omega_j(k)=(\sigma_jk)^2$ is unchanged under the transformation $k\to -k$.  Similarly, $g_\ell^{(j)}(\nu^2,T)$ and $h_\ell^{(j)}(\nu^2,T)$ are unchanged under $\nu\to-\nu$.  Hence, we have a second set of global relations\begin{equation}\label{GR_common_minus}
\begin{split}
e^{\nu^2 T}\hat{u}^{(j)}\left(\frac{-\nu}{\sigma_j},T\right)=&\hat{u}^{(j)}_0\left(\frac{-\nu}{\sigma_j}\right)+e^{\frac{i\nu x_j}{\sigma_j}} \left(\sigma_j^2 h_1^{(j)}(\nu^2,T)-i \sigma_j\nu h_0^{(j)}(\nu^2,T)\right)\\
&-e^{\frac{i\nu x_{j-1}}{\sigma_j}} \left(\sigma_j^2 g_1^{(j)}(\nu^2,T)-i \sigma_j\nu g_0^{(j)}(\nu^2,T)\right),
\end{split}
\end{equation}
where $1\leq j\leq n+1$.  In contrast to equations on an unbounded spatial domain, Equations~\eqref{GR_common_plus} and~\eqref{GR_common_minus} are valid for all $k\in\CC$.  

\begin{figure}[htbp]
\begin{center}
\def\svgwidth{.5\textwidth}
   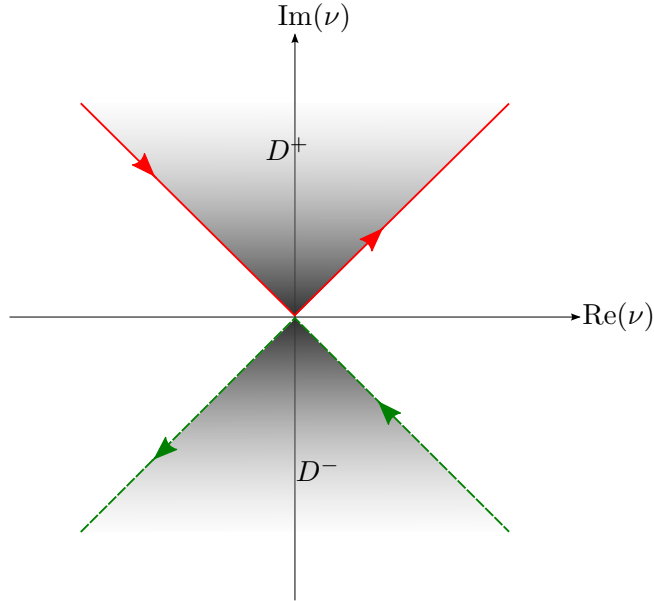 
   \caption{The regions $D^\pm$ for the heat equation.   \label{fig:heat_Dpm}}
  \end{center}
\end{figure}

\begin{figure}[htbp]
\begin{center}
\def\svgwidth{.5\textwidth}
   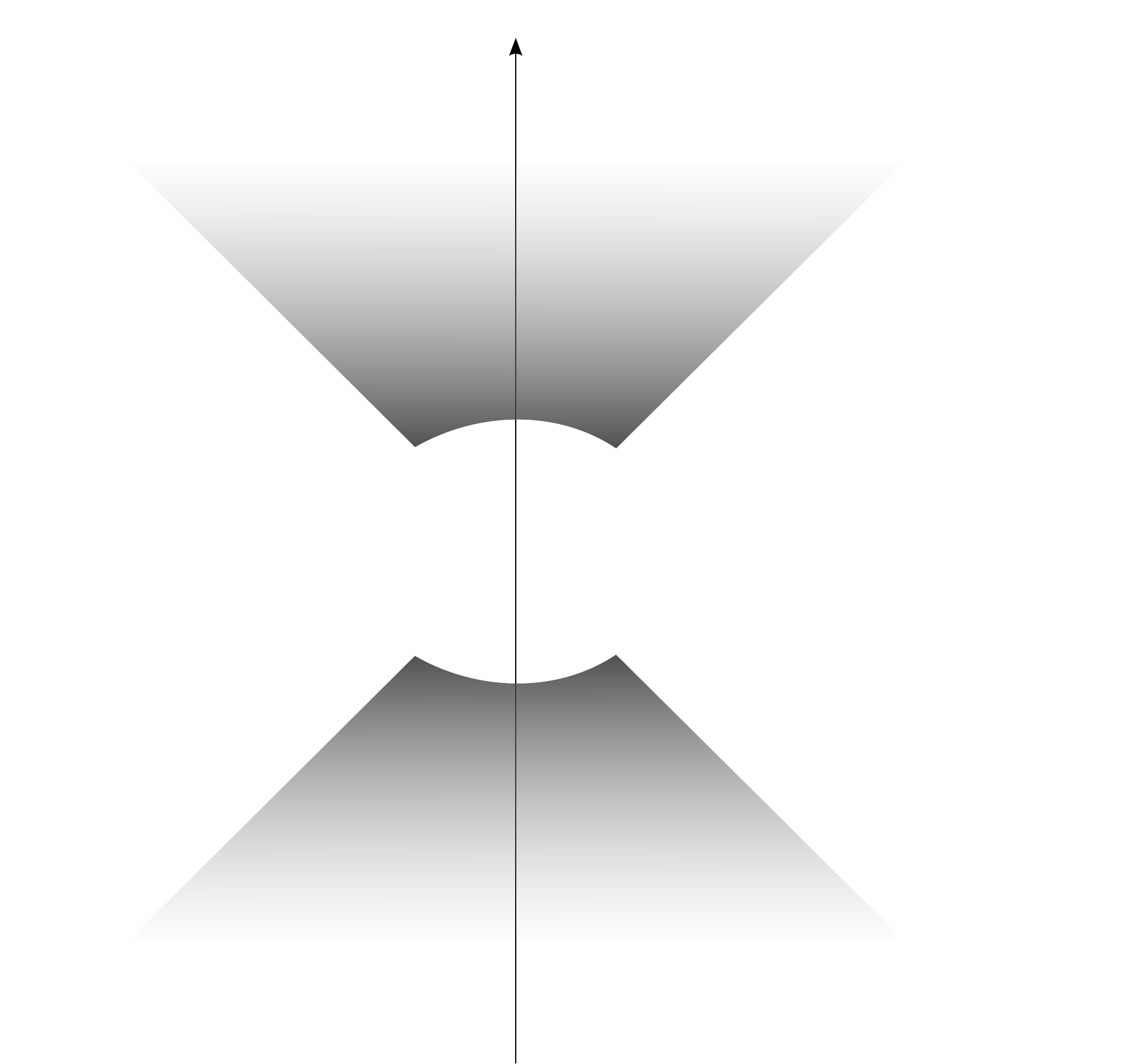 
   \caption{The regions $D_R^\pm$ for the heat equation.   \label{fig:heat_DRpm}}
  \end{center}
\end{figure}

Evaluating at $T=t$ and inverting the Fourier transform in~\eqref{GR} we have the solution formulae
\begin{equation}
\begin{split}
u^{(j)}(x,t)=&\frac{1}{2\pi}\int_{-\infty}^\infty e^{ikx-\omega_j(k) t}\hat{u}^{(j)}_0(k)\ud k\\
&+\frac{\sigma_j^2 }{2\pi}\int_{-\infty}^\infty e^{ik(x-x_j)-\omega_j(k) t} \left(h_1^{(j)}(\omega_j(k),t)+ i k h_0^{(j)}(\omega_j(k),t)\right)\ud k\\
&-\frac{\sigma_j^2 }{2\pi}\int_{-\infty}^\infty  e^{ik(x-x_{j-1})-\omega_j(k) t} \left(g_1^{(j)}(\omega_j(k),t)+ i k g_0^{(j)}(\omega_j(k),t)\right)\ud k,
\end{split}
\end{equation}
where $1\leq j\leq n+1$.  Using the change of variables $k=\nu/\sigma_j$ and replacing $t$ by $T$ in the arguments of $g_j$ and $h_j$ by noting that this is equivalent to replacing the integral $\int_0^t e^{ik^3} \frac{\partial^n}{\partial x^n} u^{(j)}(x_j,s)\ud s$ with $\int_0^T e^{ik^3} \frac{\partial^n}{\partial x^n} u^{(j)}(x_j,s)\ud s-\int_t^T e^{ik^3} \frac{\partial^n}{\partial x^n} u^{(j)}(x_j,s)\ud s$.  Using analyticity properties of the integrand and Jordan's Lemma, the contribution from the second integral is zero and thus,
\begin{equation}\label{badsoln}
\begin{split}
u^{(j)}(x,t)=&\frac{1}{2\pi}\int_{-\infty}^\infty e^{ikx-\omega_j(k) t}\hat{u}^{(j)}_0(k)\ud k\\
&+\frac{1}{2\pi}\int_{-\infty}^\infty e^{\frac{i\nu(x-x_j)}{\sigma_j}-\nu^2 t} \left(\sigma_j h_1^{(j)}(\nu^2,T)+ i \nu h_0^{(j)}(\nu^2,T)\right)\ud \nu\\
&-\frac{1}{2\pi}\int_{-\infty}^\infty  e^{\frac{i\nu(x-x_{j-1})}{\sigma_j}-\nu^2 t} \left(\sigma_j g_1^{(j)}(\nu^2,T)+ i \nu g_0^{(j)}(\nu^2,T)\right)\ud \nu,
\end{split}
\end{equation} with $1\leq j\leq n+1$ and $x_{j-1}<x<x_j$.  The integrand of the second integral in~\eqref{badsoln} is analytic and decays as $|\nu|\to\infty$ from within the set bounded between $\mathbb{R}$ and $\partial D^{-}$, and the integrand of the second integral in~\eqref{badsoln} is analytic and decays as $|\nu|\to\infty$ from within the set bounded between $\mathbb{R}$ and $\partial D^+$. Hence, by Jordan's Lemma and Cauchy's Theorem, the contours of integration can be deformed from $\mathbb{R}$ to $-\partial D^{-}$ and $\partial D^+$ respectively. \label{JordanCauchyArgument}
\begin{equation}\label{badsoln_d}
\begin{split}
u^{(j)}(x,t)=&\frac{1}{2\pi}\int_{-\infty}^\infty e^{ikx-\omega_j(k) t}\hat{u}^{(j)}_0(k)\ud k\\
&-\frac{1}{2\pi}\int_{\partial D^-} e^{\frac{i\nu(x-x_j)}{\sigma_j}-\nu^2 t} \left(\sigma_j h_1^{(j)}(\nu^2,T)+ i \nu h_0^{(j)}(\nu^2,T)\right)\ud \nu\\
&-\frac{1}{2\pi}\int_{\partial D^+}  e^{\frac{i\nu(x-x_{j-1})}{\sigma_j}-\nu^2 t} \left(\sigma_j g_1^{(j)}(\nu^2,T)+ i \nu g_0^{(j)}(\nu^2,T)\right)\ud \nu,
\end{split}
\end{equation} with $1\leq j\leq n+1$ and $x_{j-1}<x<x_j$ where $D^\pm=\{\nu\in\CC^\pm: \Re(\nu^2)<0\}$ as in Figure~\ref{fig:heat_Dpm}.  This solution is ineffective because it depends on the value of the function and its derivative evaluated at \emph{all} the interfaces and boundaries.  In order to avoid difficulties in formulas which follow, we further deform $D\pm$ to $D_R^\pm=\{\nu\in D^\pm: |\nu|>R\}$ as in Figure~\ref{fig:heat_DRpm},


\begin{equation}\label{badsoln_dR_T}
\begin{split}
u^{(j)}(x,t)=&\frac{1}{2\pi}\int_{-\infty}^\infty e^{ikx-\omega_j(k) t}\hat{u}^{(j)}_0(k)\ud k\\
&-\frac{1}{2\pi}\int_{\partial D_R^-} e^{\frac{i\nu(x-x_j)}{\sigma_j}-\nu^2 t} \left(\sigma_j h_1^{(j)}(\nu^2,T)+ i \nu h_0^{(j)}(\nu^2,T)\right)\ud \nu\\
&-\frac{1}{2\pi}\int_{\partial D_R^+}  e^{\frac{i\nu(x-x_{j-1})}{\sigma_j}-\nu^2 t} \left(\sigma_j g_1^{(j)}(\nu^2,T)+ i \nu g_0^{(j)}(\nu^2,T)\right)\ud \nu,
\end{split}
\end{equation}
for $1\leq j\leq n+1$. We replace $t$ by $T$ in the arguments of $g_j$ and $h_j$ as before \begin{equation}\label{badsoln_dR_t}
\begin{split}
u^{(j)}(x,t)=&\frac{1}{2\pi}\int_{-\infty}^\infty e^{ikx-\omega_j(k) t}\hat{u}^{(j)}_0(k)\ud k\\
&-\frac{1}{2\pi}\int_{\partial D_R^-} e^{\frac{i\nu(x-x_j)}{\sigma_j}-\nu^2 t} \left(\sigma_j h_1^{(j)}(\nu^2,T)+ i \nu h_0^{(j)}(\nu^2,T)\right)\ud \nu\\
&-\frac{1}{2\pi}\int_{\partial D_R^+}  e^{\frac{i\nu(x-x_{j-1})}{\sigma_j}-\nu^2 t} \left(\sigma_j g_1^{(j)}(\nu^2,T)+ i \nu g_0^{(j)}(\nu^2,T)\right)\ud \nu,
\end{split}
\end{equation}

While Equation~\eqref{badsoln_dR_t} makes the functional dependence of the solution more complicated than in Equation~\eqref{badsoln_dR_T}, it is useful for doing long time asymptotics, \emph{i.e.} taking the limit as $t\to\infty$.  Equation~\eqref{badsoln_dR_T} is useful for checking that the expression satisfies the equation.  While the integrands of these expressions are different, the integrals are equal and thus one may switch between them whenever convenient.

\subsection{Imperfect thermal contact}\label{sec:imperfect}
Multiplying the boundary and interface conditions~\eqref{boundary1},~\eqref{boundary2},~\eqref{ip_jump1},~\eqref{ip_jump2} by $e^{\nu^2 t}$ and integrating the result from $0$ to $T$ with respect to $t$ gives
\begin{subequations}\label{t_conds}
\begin{align}
\beta_1 g_0^{(1)}(\nu^2,T)+\beta_2 g_1^{(1)}(\nu^2,T)=\tilde{f}_1(\nu^2,T),\label{t_boundary1}\\
\beta_3 h_0^{(n+1)}(\nu^2,T)+\beta_4 h_1^{(n+1)}(\nu^2,T)=\tilde{f}_2(\nu^2,T), \label{t_boundary2}\\
\sigma_j^2 h_1^{(j)}(\nu^2,T)=H_j\left(g_0^{(j+1)}(\nu^2,T)-h_0^{(j)}(\nu^2,T)\right), \label{t_jump1}\\
\sigma_{j+1}^2 g_1^{(j+1)}(\nu^2,T)=H_j\left(g_0^{(j+1)}(\nu^2,T)-h_0^{(j)}(\nu^2,T)\right),  \label{t_jump2}
\end{align}
\end{subequations}
for $1\leq j\leq n$ where
\begin{align*}
\tilde{f}_1(\omega,t)&=\int_0^t e^{\omega s} f_1(s)\ud s,\\
\tilde{f}_2(\omega,t)&=\int_0^t e^{\omega s}f_2(s)\ud s.
\end{align*}
Applying~\eqref{t_jump1} and \eqref{t_jump2} in~\eqref{badsoln_dR_T} we have

\begin{subequations}\label{badsoln2}
\begin{equation}\label{badsoln2_1}
\begin{split}
u^{(1)}(x,t)=&\frac{1}{2\pi}\int_{-\infty}^\infty e^{ikx-\omega_1 t}\hat{u}^{(1)}_0(k)\ud k\\
&-\frac{1}{2\pi\sigma_1}\int_{\partial D_R^-} e^{\frac{i\nu(x-x_1)}{\sigma_1}-\nu^2 t} \left(H_1 g_0^{(2)}(\nu^2,T)+\left(i\sigma_1\nu-H_1 \right)h_0^{(1)}(\nu^2,T)\right)\ud \nu\\
&-\frac{1}{2\pi}\int_{\partial D_R^+}  e^{\frac{i\nu x}{\sigma_1}-\nu^2 t} \left(\sigma_1 g_1^{(1)}(\nu^2,T)+ i \nu g_0^{(1)}(\nu^2,T)\right)\ud \nu,
\end{split}
\end{equation}
\begin{equation}\label{badsoln2_j}
\begin{split}
u^{(j)}(x,t)=&\frac{1}{2\pi}\int_{-\infty}^\infty e^{ikx-\omega_j(k) t}\hat{u}^{(j)}_0(k)\ud k\\
&-\frac{1}{2\pi\sigma_j}\int_{\partial D_R^-} e^{\frac{i\nu(x-x_j)}{\sigma_j}-\nu^2 t} \left(H_j g_0^{(j+1)}(\nu^2,T)+\left(i\sigma_j\nu-H_j\right)h_0^{(j)}(\nu^2,T)\right)\ud \nu\\
&-\frac{1}{2\pi\sigma_j}\int_{\partial D_R^+}  e^{\frac{i\nu(x-x_{j-1})}{\sigma_j}-\nu^2 t} \left(\left(H_{j-1}+ i \sigma_j \nu\right) g_0^{(j)}(\nu^2,T)-H_{j-1}h_0^{(j-1)}(\nu^2,T)\right)\ud \nu,
\end{split}
\end{equation}
\begin{equation}\label{badsoln2_n1}
\begin{split}
u^{(n+1)}(x,t)=&\frac{1}{2\pi}\int_{-\infty}^\infty e^{ikx-\omega_{n+1} t}\hat{u}^{(n+1)}_0(k)\ud k\\
&-\frac{1}{2\pi}\int_{\partial D_R^-} e^{\frac{i\nu(x-x_{n+1})}{\sigma_{n+1}}-\nu^2 t} \left(\sigma_{n+1} h_1^{(n+1)}(\nu^2,T)+ i \nu h_0^{(n+1)}(\nu^2,T)\right)\ud \nu\\
&-\frac{1}{2\pi\sigma_{n+1}}\int_{\partial D_R^+}  e^{\frac{i\nu(x-x_{n})}{\sigma_{n+1}}-\nu^2 t} \left(\left(H_{n}+i\sigma_{n+1}\nu\right)g_0^{(n+1)}(\nu^2,T)-H_{n}h_0^{(n)}(\nu^2,T)\right)\ud \nu,
\end{split}
\end{equation}
\end{subequations}
 with $2\leq j\leq n$.

We use~\eqref{t_jump1} and \eqref{t_jump2} in the global relations~\eqref{GR_common_plus} and~\eqref{GR_common_minus} which gives
\begin{subequations}\label{new_GRs}
\begin{equation}
\begin{split}
e^{\nu^2 T}\hat{u}^{(1)}\left(\frac{\nu}{\sigma_1},T\right)=&\hat{u}^{(1)}_0\left(\frac{\nu}{\sigma_1}\right)+e^{-\frac{i\nu x_1}{\sigma_1}} \left(H_1g_0^{(2)}(\nu^2,T)+( i \sigma_1\nu -H_1) h_0^{(1)}(\nu^2,T)\right)\\
&-\sigma_1^2 g_1^{(1)}(\nu^2,T)- i \sigma_1\nu g_0^{(1)}(\nu^2,T),
\end{split}
\end{equation}
\begin{equation}
\begin{split}
e^{\nu^2 T}\hat{u}^{(j)}\left(\frac{\nu}{\sigma_j},T\right)=&\hat{u}^{(j)}_0\left(\frac{\nu}{\sigma_j}\right)+e^{-\frac{i\nu x_j}{\sigma_j}} \left(H_j g_0^{(j+1)}(\nu^2,T)+( i \sigma_j\nu -H_j) h_0^{(j)}(\nu^2,T)\right)\\
&+e^{-\frac{i\nu x_{j-1}}{\sigma_j}} \left(H_{j-1}h_0^{(j-1)}(\nu^2,T)-(H_{j-1}+ i \sigma_j\nu)g_0^{(j)}(\nu^2,T)\right),
\end{split}
\end{equation}
\begin{equation}
\begin{split}
e^{\nu^2 T}\hat{u}^{(n+1)}\left(\frac{\nu}{\sigma_{n+1}},T\right)=&\hat{u}^{(n+1)}_0\left(\frac{\nu}{\sigma_{n+1}}\right)+e^{-\frac{i\nu x_{n+1}}{\sigma_{n+1}}} \left(\sigma_{n+1}^2 h_1^{(n+1)}(\nu^2,T)+ i \sigma_{n+1}\nu h_0^{(n+1)}(\nu^2,T)\right)\\
&+e^{-\frac{i\nu x_{n}}{\sigma_{n+1}}} \left(H_{n}h_0^{(n)}(\nu^2,T)-(H_{n}+ i \sigma_{n+1}\nu)g_0^{(n+1)}(\nu^2,T)\right),
\end{split}
\end{equation}
\begin{equation}
\begin{split}
e^{\nu^2 T}\hat{u}^{(1)}\left(\frac{-\nu}{\sigma_1},T\right)=&\hat{u}^{(1)}_0\left(\frac{-\nu}{\sigma_1}\right)+e^{\frac{i\nu x_1}{\sigma_1}} \left(H_1g_0^{(2)}(\nu^2,T)-(i \sigma_1\nu+H_1) h_0^{(1)}(\nu^2,T)\right)\\
&- \sigma_1^2 g_1^{(1)}(\nu^2,T)+i \sigma_1\nu g_0^{(1)}(\nu^2,T),
\end{split}
\end{equation}
\begin{equation}
\begin{split}
e^{\nu^2 T}\hat{u}^{(j)}\left(\frac{-\nu}{\sigma_j},T\right)=&\hat{u}^{(j)}_0\left(\frac{-\nu}{\sigma_j}\right)+e^{\frac{i\nu x_j}{\sigma_j}} \left(H_jg_0^{(j+1)}(\nu^2,T)-(i \sigma_j\nu+H_{j}) h_0^{(j)}(\nu^2,T)\right)\\
&+e^{\frac{i\nu x_{j-1}}{\sigma_j}} \left(H_{j-1}h_0^{(j-1)}(\nu^2,T)-(H_{j-1}-i\sigma_j\nu)g_0^{(j)}(\nu^2,T)\right),
\end{split}
\end{equation}
\begin{equation}
\begin{split}
e^{\nu^2 T}\hat{u}^{(n+1)}\left(\frac{-\nu}{\sigma_{n+1}},T\right)=&\hat{u}^{(n+1)}_0\left(\frac{-\nu}{\sigma_{n+1}}\right)+e^{\frac{i\nu x_{n+1}}{\sigma_{n+1}}} \left(\sigma_{n+1}^2 h_1^{(n+1)}(\nu^2,T)-i \sigma_{n+1}\nu h_0^{(n+1)}(\nu^2,T)\right)\\
&+e^{\frac{i\nu x_{n}}{\sigma_{n+1}}} \left(H_{n}h_0^{(n)}(\nu^2,T)-(H_{n}-i\sigma_{n+1}\nu)g_0^{(n+1)}(\nu^2,T)\right),
\end{split}
\end{equation}
\end{subequations}
where $2\leq j\leq n$.  Equation~\eqref{badsoln2} involves $2n+4$ unknown functions $g_0^{(j)}(\nu^2,T)$, $h_0^{(j)}(\nu^2,T)$, $g_1^{(1)}(\nu^2,T)$, $h_1^{(n)}(\nu^2,t)$ for $1\leq j\leq n+1.$  However, these same unknown functions are related through the $2n+2$ global relations~\eqref{new_GRs} and the transformed boundary conditions~\eqref{t_boundary1} and~\eqref{t_boundary2}.  Solving this linear system for the unknown functions amounts to solving the $(2n+4)\times (2n+4)$ matrix problem $$\mathcal{A}(\nu)X(\nu^2,T)=\mathcal{Y}(\nu,T)+Y(\nu,T)$$ where

\begin{subequations}
{\footnotesize
\begin{equation}\label{X}
X(\nu^2,T)=\left(g_1^{(1)}(\nu^2,T), g_0^{(1)}(\nu^2,T),\ldots,g_0^{(n+1)}(\nu^2,T),h_0^{(1)}(\nu^2,T),\ldots,h_0^{(n+1)}(\nu^2,T), h_1^{(n+1)}(\nu^2,T)\right)^\top,
\end{equation}
\begin{equation}\label{Y}
Y(\nu,T)=-\left(-\tilde{f}_1(\nu^2,T),\hat{u}_0^{(1)}\left(\frac{\nu}{\sigma_1}\right),\cdots,\hat{u}_0^{(n+1)}\left(\frac{\nu}{\sigma_{n+1}}\right),\hat{u}_0^{(1)}\left(\frac{-\nu}{\sigma_1}\right),\cdots,\hat{u}_0^{(n+1)}\left(\frac{-\nu}{\sigma_{n+1}}\right), -\tilde{f}_2(\nu^2,T)\right)^\top,
\end{equation}
\begin{equation}\label{badY}
\mathcal{Y}(\nu,T)=e^{\nu^2 T}\left(0,\hat{u}^{(1)}\left(\frac{\nu}{\sigma_1},T\right),\cdots,\hat{u}^{(n+1)}\left(\frac{\nu}{\sigma_{n+1}},T\right),\hat{u}^{(1)}\left(\frac{-\nu}{\sigma_1},T\right),\cdots,\hat{u}^{(n+1)}\left(\frac{-\nu}{\sigma_{n+1}},T\right), 0\right)^\top,
\end{equation}
}
and 
{\footnotesize
$$\mathcal{A}_{11}(\nu)=\left(\begin{array}{cccccc}
\beta_2&\beta_1\\
-\sigma_1^2&-i\sigma_1\nu &H_1e^{-i\frac{\nu x_1}{\sigma_1}}\\
0&0&-(H_1+i\sigma_2\nu) e^{-i\frac{\nu x_1}{\sigma_2}}&H_2e^{-i\frac{\nu x_2}{\sigma_2}}\\
\vdots&&\hspace{1.3in}\ddots&\hspace{1.5in}\ddots&\\
\vdots&&&-(H_{n-1}+i\sigma_{n}\nu) e^{-i\frac{\nu x_{n-1}}{\sigma_{n}}}&H_{n}e^{-i\frac{\nu x_{n}}{\sigma_{n}}}\\
0&\hdots&\hdots&0&-(H_{n}+i\sigma_{n+1}\nu) e^{-i\frac{\nu x_{n}}{\sigma_{n+1}}}\\

\end{array}\right),
$$
$$\mathcal{A}_{12}(\nu)=\left(\begin{array}{cccccc}
0&\hdots&\hdots&0&0\\
(i\sigma_1\nu-H_1)e^{-i\frac{\nu x_1}{\sigma_1}}&&&&0\\
H_1e^{-i\frac{\nu x_1}{\sigma_2}}&(i\sigma_2\nu-H_2)e^{-i\frac{\nu x_2}{\sigma_2}}&0&\cdots&0\\
\hspace{.5in}\ddots&\hspace{.5in}\ddots&&&\vdots\\
&H_{n-1}e^{-i\frac{\nu x_{n-1}}{\sigma_{n}}}&(i\sigma_{n}\nu-H_{n})e^{-i\frac{\nu x_{n}}{\sigma_{n}}}&&\vdots\\
&&H_{n}e^{-i\frac{\nu x_{n}}{\sigma_{n+1}}}&i\sigma_{n+1}\nu e^{-i\frac{\nu x_{n+1}}{\sigma_{n+1}}} &\sigma_{n+1}^2e^{-i\frac{\nu x_{n+1}}{\sigma_{n+1}}}

\end{array}\right),
$$

$$\mathcal{A}_{21}(\nu)=\left(\begin{array}{cccccc}
-\sigma_1^2&i\sigma_1\nu &H_1e^{i\frac{\nu x_1}{\sigma_1}}&0\\
0&0&(i\sigma_2\nu-H_1) e^{i\frac{\nu x_1}{\sigma_2}}&H_2e^{i\frac{\nu x_2}{\sigma_2}}\\
\vdots&&\hspace{1.3in}\ddots&\hspace{1.5in}\ddots&\\
\vdots&&&(i\sigma_{n}\nu-H_{n-1}) e^{i\frac{\nu x_{n-1}}{\sigma_{n}}}&H_{n}e^{i\frac{\nu x_{n}}{\sigma_{n}}}\\
0&&&&(i\sigma_{n+1}\nu-H_{n}) e^{i\frac{\nu x_{n}}{\sigma_{n+1}}}\\
0&\hdots&\hdots&\hdots&0
\end{array}\right),
$$
$$\mathcal{A}_{22}(\nu)=\left(\begin{array}{cccccc}
-(i\sigma_1\nu+H_1)e^{i\frac{\nu x_1}{\sigma_1}}\\
H_1e^{i\frac{\nu x_1}{\sigma_2}}&-(i\sigma_2\nu+H_2)e^{i\frac{\nu x_2}{\sigma_2}}\\
\hspace{.5in}\ddots&\hspace{.5in}\ddots&\\
&H_{n-1}e^{i\frac{\nu x_{n-1}}{\sigma_{n}}}&-(i\sigma_{n}\nu+H_{n})e^{i\frac{\nu x_{n}}{\sigma_{n}}}\\
&&H_{n}e^{i\frac{\nu x_{n}}{\sigma_{n+1}}}&-i\sigma_{n+1}\nu e^{i\frac{\nu x_{n+1}}{\sigma_{n+1}}} &\sigma_{n+1}^2e^{i\frac{\nu x_{n+1}}{\sigma_{n+1}}}\\
&&&\beta_3&\beta_4\\
\end{array}\right),
$$
}
and 
\begin{equation}\label{Amatrix}
\mathcal{A}(\nu)=\left(\begin{array}{c:c}
\mathcal{A}_{11}(\nu)&\mathcal{A}_{12}(\nu)\\
\hdashline
\mathcal{A}_{21}(\nu)&\mathcal{A}_{22}(\nu)
\end{array}
\right).
\end{equation}
\end{subequations}

The matrix $\mathcal{A}(\nu)$ is made up of four $(n+2)\times(n+2)$ blocks as indicated by the dashed lines.  $\mathcal{A}_{11}$ has nonzero entries only on the main and $+1$ diagonals, $\mathcal{A}_{12}$ has nonzero entries on the $-1$ and $-2$ diagonals, $\mathcal{A}_{21}$ has nonzero entries on the $+1$ and $+2$ diagonals, and $\mathcal{A}_{22}$ has nonzero entries on the main and $-1$ diagonals.  The boundary conditions are incorporated in the first and last rows of $\mathcal{A}(\nu)$.


The matrix $\mathcal{A}(\nu)$ is singular for isolated values of $\nu$.  Asymptotically, for large $|\nu|$, the zeros of $\det(\mathcal{A}(\nu))$ are on the real line~\cite{Langer}.  Since asymptotically there are no zeros in $D_R^+$, a sufficiently large $R$ may be chosen such that $\mathcal{A}(\nu)$ is nonsingular for every $\nu\in D_R^+$ and $\det(\mathcal{A}(\nu))\neq0$. 

Every term in the linear equation $\mathcal{A}(\nu)X(\nu^2,T)=Y(\nu,T)$ is known.  By substituting the solutions of this equation into~\eqref{badsoln2}, we have solved the heat equation on the finite interval with $n$ interfaces with imperfect interface conditions in terms of only known functions.  It remains to show that the contribution to the solution from the linear equation $\mathcal{A}(\nu)X(\nu^2,T)=\mathcal{Y}(\nu,T)$ is 0 when substituted into~\eqref{badsoln2}.

To this end consider $\mathcal{A}(\nu)X(\nu^2,T)=\mathcal{Y}(\nu,T).$  For the integral over $\partial D_R^+$ we factor $\mathcal{A}(\nu)=\mathcal{A}^{(L,+)}(\nu)\mathcal{A}^{(M,+)}(\nu)$ where
$$\mathcal{A}^{(L,+)}(\nu)=\left(\begin{array}{cccc:cccc}
1\\
&e^{-i\frac{\nu x_1}{\sigma_1}}&&&\\
&&\ddots\\
&&&e^{-i\frac{\nu x_{n+1}}{\sigma_{n+1}}}&\\
\hdashline
&&&&e^{i\frac{\nu x_0}{\sigma_1}}\\
&&&&&\ddots\\
&&&&&&e^{i\frac{\nu x_{n}}{\sigma_{n+1}}}\\
&&&&&&&1
\end{array}\right).
$$
For the integral over $\partial D_R^-$ we factor $\mathcal{A}(\nu)=\mathcal{A}^{(L,-)}(\nu)\mathcal{A}^{(M,-)}(\nu)$,
where
$$\mathcal{A}^{(L,-)}(\nu)=\left(\begin{array}{cccc:cccc}
1\\
&e^{i\frac{\nu x_1}{\sigma_1}}&&&\\
&&\ddots\\
&&&e^{i\frac{\nu x_{n+1}}{\sigma_{n+1}}}&\\
\hdashline
&&&&e^{i\frac{\nu x_1}{\sigma_1}}\\
&&&&&\ddots\\
&&&&&&e^{i\frac{\nu x_{n+1}}{\sigma_{n+1}}}\\
&&&&&&&1
\end{array}\right).
$$
Let $\mathcal{A}_j(\nu,T)$ be the matrix $\mathcal{A}(\nu)$ with the $j^{\textrm{th}}$ column replaced by $\mathcal{Y}(\nu,T)$.  Similar to $\mathcal{A}(\nu )$, this matrix can be factored as $\mathcal{A}_j(\nu,T)=\mathcal{A}^{(L,\pm)}(\nu )\mathcal{A}_j^{(M,\pm)}(\nu,T)\mathcal{A}_j^R(\nu,T)$ where $\mathcal{A}_j^R(\nu,T)$ is the $(2n+4)\times (2n+4)$ identity matrix with the $(j,j)$ entry replaced by $e^{\nu^2 T}$.  Hence, $\det(\mathcal{A}_j(\nu,T))=e^{\nu ^2T}\det(\mathcal{A}^{(L,\pm)}(\nu))\det(\mathcal{A}_j^{(M,\pm)}(\nu,T))$.

The terms we are trying to eliminate contribute to the solution~\eqref{badsoln2} in the form:
\begin{subequations}
\begin{align}
&-\frac{1}{2\pi\sigma_j}\int_{\partial D_R^-} e^{\frac{i\nu(x-x_j)}{\sigma_j}-\nu^2 t} \left(H_j g_0^{(j+1)}(\nu^2,T)+\left(i\sigma_j\nu-H_j\right)h_0^{(j)}(\nu^2,T)\right)\ud \nu,\\
&-\frac{1}{2\pi}\int_{\partial D_R^-} e^{\frac{i\nu(x-x_n)}{\sigma_n}-\nu^2 t} \left(\sigma_n h_1^{(n)}(\nu^2,T)+ i \nu h_0^{(n)}(\nu^2,T)\right)\ud \nu,
\end{align}
\end{subequations}
for $1\leq j\leq n-1$, and
\begin{subequations}
\begin{align}
&-\frac{1}{2\pi}\int_{\partial D_R^+}  e^{\frac{i\nu x}{\sigma_1}-\nu^2 t} \left(\sigma_1  g_1^{(1)}(\nu^2,T)+ i \nu g_0^{(1)}(\nu^2,T)\right)\ud \nu,\\
&-\frac{1}{2\pi\sigma_j}\int_{\partial D_R^+}  e^{\frac{i\nu(x-x_{j-1})}{\sigma_j}-\nu^2 t} \left(\left(H_{j-1}+ i \sigma_j \nu\right) g_0^{(j)}(\nu^2,T)-H_{j-1}h_0^{(j-1)}(\nu^2,T)\right)\ud \nu,
\end{align}
\end{subequations}
for $2\leq j\leq n$ with $x_{j-1}<x<x_j$.  Using Cramer's Rule these become
\begin{subequations}\label{badpart_minus}
\begin{align}
&-\frac{1}{2\pi\sigma_j}\int_{\partial D_R^-} e^{\frac{i\nu(x-x_j)}{\sigma_j}+\nu^2(T- t)} \left(H_j \frac{\det(\mathcal{A}_{j+2}^{(M,-)})}{\det(\mathcal{A}^{(M,-)})}+\left(i\sigma_j\nu-H_j\right)\frac{\det(\mathcal{A}_{n+j+1}^{(M,-)})}{\det(\mathcal{A}^{(M,-)})} \right)\ud \nu,\label{badpart_minus_a}\\ \label{badpart_minus2}
&-\frac{1}{2\pi}\int_{\partial D_R^-} e^{\frac{i\nu(x-x_n)}{\sigma_n}+\nu^2 (T-t)} \left(\sigma_n \frac{\det(\mathcal{A}_{2n+2}^{(M,-)})}{\det(\mathcal{A}^{(M,-)})}+ i \nu \frac{\det(\mathcal{A}_{2n+1}^{(M,-)})}{\det(\mathcal{A}^{(M,-)})}\right)\ud \nu,
\end{align}
\end{subequations}
for $1\leq j\leq n-1$, and
\begin{subequations}\label{badpart_plus}
\begin{align}
&-\frac{1}{2\pi}\int_{\partial D_R^+}  e^{\frac{i\nu x}{\sigma_1}+\nu^2 (T-t)} \left(\sigma_1 \frac{\det(\mathcal{A}_{1}^{(M,+)})}{\det(\mathcal{A}^{(M,+)})}+ i \nu \frac{\det(\mathcal{A}_{2}^{(M,+)})}{\det(\mathcal{A}^{(M,+)})}  \right)\ud \nu,\\ \label{badpart_plus2}
&-\frac{1}{2\pi\sigma_j}\int_{\partial D_R^+}  e^{\frac{i\nu(x-x_{j-1})}{\sigma_j}+\nu^2(T- t)} \left(\left(H_{j-1}+ i \sigma_j \nu\right) \frac{\det(\mathcal{A}_{j+1}^{(M,+)})}{\det(\mathcal{A}^{(M,+)})} -H_{j-1} \frac{\det(\mathcal{A}_{n+j}^{(M,+)})}{\det(\mathcal{A}^{(M,+)})}  \right)\ud \nu,
\end{align}
\end{subequations}
for $2\leq j\leq n$.  As usual in the Fokas Method we use the large $\nu$ asymptotics to show the terms in~\eqref{badpart_minus} and~\eqref{badpart_plus} are 0.  Observe the elements of $\mathcal{A}^{(M,\pm)}$ are either 0, $\mathcal{O}(\nu)$, or decaying exponentially fast for $\nu\in D^\pm$ respectively.  Hence, $$\det(\mathcal{A}^{(M,\pm)}(\nu))=c(\nu)=\mathcal{O}(\nu^{2n+2})$$ for large $\nu$ in $D_R^\pm$.

\begin{figure}[htbp]
\begin{center}
\def\svgwidth{.5\textwidth}
   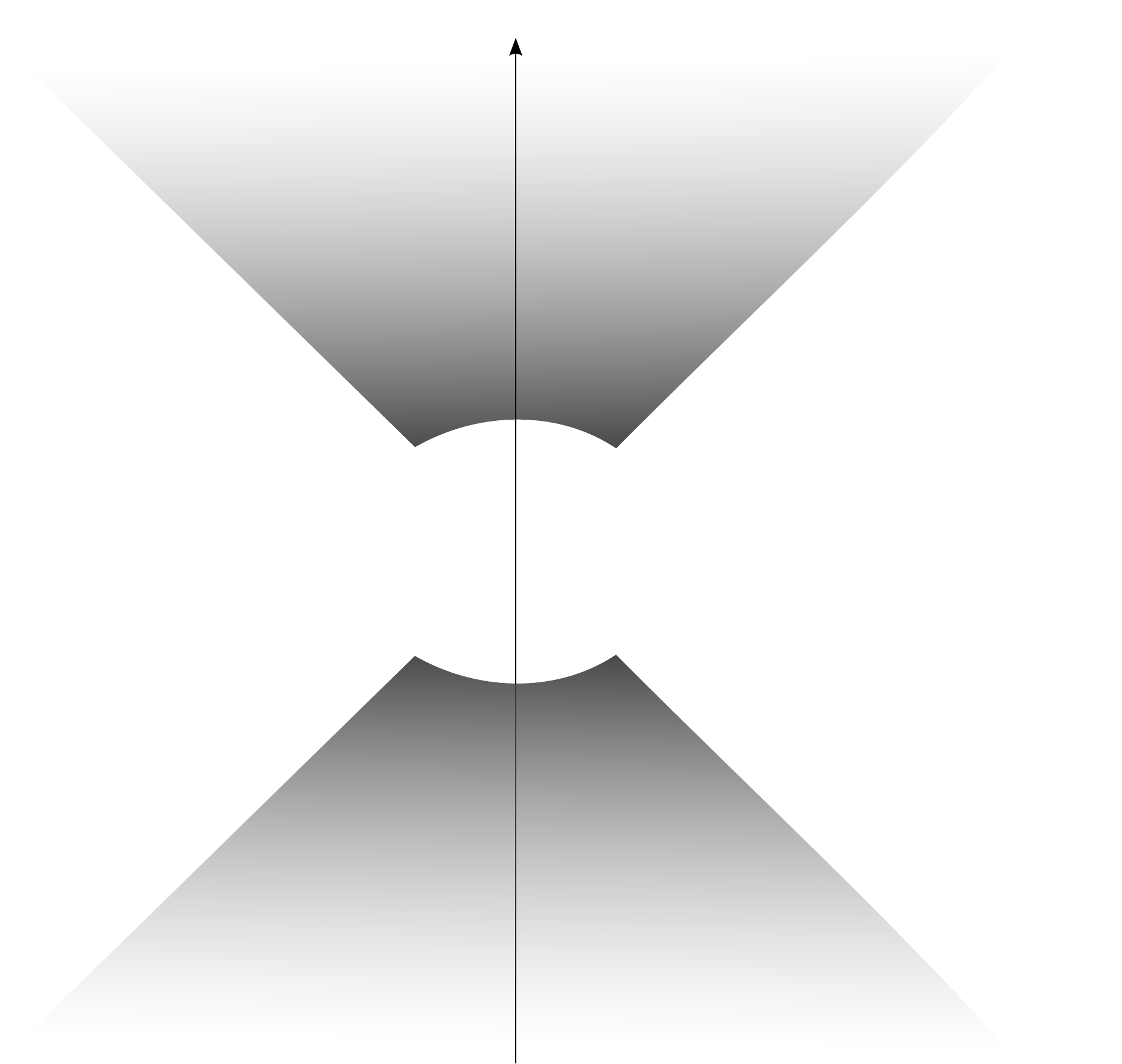 
   \caption{  \label{fig:heat_DRpm_close}}
  \end{center}
\end{figure}

We begin by examining the first term of~\eqref{badpart_minus_a}. Expanding the determinant of $\mathcal{A}^{(M,-)}_j(\nu,t)$ along the $j^{\textrm{th}}$ column we see that 

\begin{align}
&e^{i\frac{\nu(x-x_j)}{\sigma_j}+\nu^2(T-t)} \frac{\det(\mathcal{A}_{j+2}^{(M,-)})}{\det(\mathcal{A}^{(M,-)})}= e^{i\frac{\nu(x-x_j)}{\sigma_j}+\nu^2(T-t)} \frac{\det(\mathcal{A}_{j+2}^{(M,-)})}{c(\nu)}\nonumber\\
&=e^{i\frac{\nu(x-x_j)}{\sigma_j}+\nu^2(T-t)} \sum_{\ell=1}^n \left(\alpha_{\ell}(\nu)e^{-i\frac{\nu x_\ell}{\sigma_\ell}}\hat{u}^{(\ell)}\left(\frac{\nu}{\sigma_\ell},T\right)+ \beta_{\ell}(\nu)e^{i\frac{\nu x_{\ell-1}}{\sigma_{\ell} }}\hat{u}^{(\ell)}\left(-\frac{\nu}{\sigma_\ell},T\right)\right),
\end{align}
where $\alpha_{\ell}(\nu)$ and $\beta_{\ell}(\nu)$ are $\mathcal{O}(\nu^{0})$ for large $\nu$ and $x_{j-1}<x<x_j$.  Note that $$e^{i\frac{\nu(x-x_j)}{\sigma_j}+\nu^2(T-t)}e^{-i\frac{\nu x_\ell}{\sigma_\ell}}\hat{u}^{(\ell)}\left(\frac{\nu}{\sigma_\ell},T\right)$$ and $$e^{i\frac{\nu(x-x_j)}{\sigma_j}+\nu^2(T-t)}e^{i\frac{\nu x_{\ell-1}}{\sigma_{\ell} }}\hat{u}^{(\ell)}\left(-\frac{\nu}{\sigma_\ell},T\right)$$ decay exponentially fast for $|\nu|\to\infty$ from within $D_R^-$.  Thus, by Jordan's Lemma, the integrals of $e^{i\frac{\nu(x-x_j)}{\sigma_j}+\nu^2(T-t)}e^{\frac{-i\nu x_\ell}{\sigma_\ell}}\hat{u}^{(\ell)}\left(\nu/\sigma_\ell,T\right)$ and $e^{i\frac{\nu(x-x_j)}{\sigma_j}+\nu^2(T-t)}e^{\frac{i\nu x_{\ell-1}}{\sigma_{\ell-1}}}\hat{u}^{(\ell)}\left(-\nu/\sigma_\ell,T\right)$ along a closed, bounded curve in the lower-half of the complex $\nu$ plane vanish for $x_{\ell-1}<x<x_\ell$. In particular we consider the closed curve  $\mathcal{L}^{-}=\mathcal{L}_{D^{-}}\cup\mathcal{L}^{-}_C$ where $\mathcal{L}_{D^{-}}=\partial D_R^{-} \cap \{\nu: |\nu|<C\}$ and $\mathcal{L}^{-}_C=\{\nu\in D_R^{-}: |\nu|=C\}$, see Figure~\ref{fig:heat_DRpm_close}

Similar to the argument on page~\pageref{JordanCauchyArgument}, since the integral along $\mathcal{L}^{+}_C$ vanishes for large $C$, the integrals~\eqref{badpart_minus_a} must vanish since the contour $\mathcal{L}_{D^{+}}$  becomes $\partial D^{+}$ as $C\to\infty$.  The uniform decay of the ratios of the determinants for large $\nu$ is exactly the condition required for the integral to vanish using Jordan's Lemma.  This argument can be repeated for~\eqref{badpart_minus2}-\eqref{badpart_plus2}  Hence, the solution to~\eqref{heatgeneral} is~\eqref{badsoln2} where $g_1^{(1)}(\nu^2,T), g_0^{(j)}(\nu^2,T), h_0^{(j)}(\nu^2,T),$ and $h_1^{(n)}(\nu^2,T)$ for $1\leq j\leq n+1$ are found by solving
\begin{equation}\label{linear_system}
\mathcal{A}(\nu)X(\nu^2,T)=Y(\nu,T),
\end{equation} 
where $\mathcal{A}(\nu), X(\nu^2,T)$, and $Y(\nu)$ are given in Equations~\eqref{Amatrix},~\eqref{X}, and~\eqref{Y} respectively.  

\subsection{Perfect thermal contact}\label{sec:perfect}
In this subsection we will repeat much of the analysis from~\ref{sec:imperfect} for different interface conditions and generalize what is presented in~\cite{DeconinckPelloniSheils} to $n$ interfaces.

Multiplying the interface conditions~\eqref{p_jump1},~\eqref{p_jump2} by $e^{\nu^2 t}$ and integrating the result from $0$ to $T$ with respect to $t$ gives
\begin{subequations}\label{p_t_conds}
\begin{align}
h_0^{(j)}(\nu^2,T)&=g_0^{(j+1)}(\nu^2,T), \label{p_t_jump1}\\
\sigma_{j}^2 h_1^{(j)}(\nu^2,T)&=\sigma_{j+1}^2 g_1^{(j+1)}(\nu^2,T),  \label{p_t_jump2}
\end{align}
\end{subequations}
for $1\leq j\leq n$.  Applying~\eqref{p_t_jump1} and \eqref{p_t_jump2} in~\eqref{badsoln_dR_T}
\begin{subequations}\label{p_badsoln2}
\begin{equation}\label{p_badsoln2_j}
\begin{split}
u^{(j)}(x,t)=&\frac{1}{2\pi}\int_{-\infty}^\infty e^{ikx-\omega_j(k) t}\hat{u}^{(j)}_0(k)\ud k\\
&-\frac{1}{2\pi}\int_{\partial D_R^-} e^{\frac{i\nu(x-x_j)}{\sigma_j}-\nu^2 t} \left(\frac{\sigma_{j+1}^2}{\sigma_j} g_1^{(j+1)}(\nu^2,T)+ i \nu g_0^{(j+1)}(\nu^2,T)\right)\ud \nu\\
&-\frac{1}{2\pi}\int_{\partial D_R^+}  e^{\frac{i\nu(x-x_{j-1})}{\sigma_j}-\nu^2 t} \left(\sigma_j g_1^{(j)}(\nu^2,T)+ i \nu g_0^{(j)}(\nu^2,T)\right)\ud \nu,
\end{split}
\end{equation}
\begin{equation}\label{p_badsoln2_n1}
\begin{split}
u^{(n+1)}(x,t)=&\frac{1}{2\pi}\int_{-\infty}^\infty e^{ikx-\omega_{n+1}(k) t}\hat{u}^{(n+1)}_0(k)\ud k\\
&-\frac{1}{2\pi}\int_{\partial D_R^-} e^{\frac{i\nu(x-x_{n+1})}{\sigma_{n+1}}-\nu^2 t} \left(\sigma_{n+1} h_1^{(n+1)}(\nu^2,T)+ i \nu h_0^{(n+1)}(\nu^2,T)\right)\ud \nu\\
&-\frac{1}{2\pi}\int_{\partial D_R^+}  e^{\frac{i\nu(x-x_{n})}{\sigma_{n+1}}-\nu^2 t} \left(\sigma_{n+1} g_1^{(n+1)}(\nu^2,T)+ i \nu g_0^{(n+1)}(\nu^2,T)\right)\ud \nu,
\end{split}
\end{equation}
\end{subequations}
with $1\leq j\leq n$.

We use~\eqref{p_t_jump1} and \eqref{p_t_jump2} in the global relations~\eqref{GR_common_plus} and~\eqref{GR_common_minus} which gives
\begin{subequations}\label{p_new_GRs}
\begin{equation}
\begin{split}
e^{\nu^2 T}\hat{u}^{(j)}\left(\frac{\nu}{\sigma_j},T\right)=&\hat{u}^{(j)}_0\left(\frac{\nu}{\sigma_j}\right)+e^{-\frac{i\nu x_j}{\sigma_j}} \left(\sigma_{j+1}^2 g_1^{(j+1)}(\nu^2,T)+ i \sigma_j\nu g_0^{(j+1)}(\nu^2,T)\right)\\
&-e^{-\frac{i\nu x_{j-1}}{\sigma_j}} \left(\sigma_j^2 g_1^{(j)}(\nu^2,T)+ i \sigma_j\nu g_0^{(j)}(\nu^2,T)\right),
\end{split}
\end{equation}
\begin{equation}
\begin{split}
e^{\nu^2 T}\hat{u}^{(j)}\left(\frac{-\nu}{\sigma_j},T\right)=&\hat{u}^{(j)}_0\left(\frac{-\nu}{\sigma_j}\right)+e^{\frac{i\nu x_j}{\sigma_j}} \left(\sigma_{j+1}^2 g_1^{(j+1)}(\nu^2,T)-i \sigma_j\nu g_0^{(j+1)}(\nu^2,T)\right)\\
&-e^{\frac{i\nu x_{j-1}}{\sigma_j}} \left(\sigma_j^2 g_1^{(j)}(\nu^2,T)-i \sigma_j\nu g_0^{(j)}(\nu^2,T)\right),
\end{split}
\end{equation}
\begin{equation}
\begin{split}
e^{\nu^2 T}\hat{u}^{(n+1)}\left(\frac{\nu}{\sigma_{n+1}},T\right)=&\hat{u}^{(n+1)}_0\left(\frac{\nu}{\sigma_{n+1}}\right)+e^{-\frac{i\nu x_{n+1}}{\sigma_{n+1}}} \left(\sigma_{n+1}^2 h_1^{(n+1)}(\nu^2,T)+ i \sigma_{}\nu h_0^{(n+1)}(\nu^2,T)\right)\\
&-e^{-\frac{i\nu x_{n}}{\sigma_{n+1}}} \left(\sigma_{n+1}^2 g_1^{(n+1)}(\nu^2,T)+ i \sigma_{n+1}\nu g_0^{(n+1)}(\nu^2,T)\right),
\end{split}
\end{equation}
\begin{equation}
\begin{split}
e^{\nu^2 T}\hat{u}^{(n+1)}\left(\frac{-\nu}{\sigma_{n+1}},T\right)=&\hat{u}^{(n+1)}_0\left(\frac{-\nu}{\sigma_{n+1}}\right)+e^{\frac{i\nu x_{n+1}}{\sigma_{n+1}}} \left(\sigma_{n+1}^2 h_1^{(n+1)}(\nu^2,T)-i \sigma_{n+1}\nu h_0^{(n+1)}(\nu^2,T)\right)\\
&-e^{\frac{i\nu x_{n}}{\sigma_{n+1}}} \left(\sigma_{n+1}^2 g_1^{(n+1)}(\nu^2,T)-i \sigma_{n+1}\nu g_0^{(n+1)}(\nu^2,T)\right),
\end{split}
\end{equation}
\end{subequations}
where $1\leq j\leq n$.  Equation~\eqref{p_badsoln2} involves $2n+4$ unknown functions $g_0^{(j)}(\nu^2,T)$, $g_1^{(j)}(\nu^2,T)$, $h_0^{(n+1)}(\nu^2,t)$, $h_1^{(n+1)}(\nu^2,t)$ for $1\leq j\leq n+1.$  These functions are related through the $2n+2$ global relations~\eqref{p_new_GRs} and the transformed boundary conditions~\eqref{t_boundary1} and~\eqref{t_boundary2}.  Solving this linear system for the unknown functions amounts to solving the $(2n+4)\times (2n+4)$ matrix problem $$\mathcal{A}^{(p)}(\nu)X^{(p)}(\nu^2,T)=\mathcal{Y}(\nu,T)+Y(\nu,T)$$ where
\begin{subequations}
{\footnotesize
\begin{equation}\label{p_X}
X^{(p)}(\nu^2,T)=\left(g_0^{(1)}(\nu^2,T), \ldots, g_0^{(n+1)}(\nu^2,T), h_0^{(n+1)}(\nu^2,T),g_1^{(1)}(\nu^2,T),\ldots, g_1^{(n+1)}(\nu^2,T),h_1^{(n+1)}(\nu^2,T)\right)^\top,
\end{equation}
}
and
{\footnotesize
$$\mathcal{A}^{(p)}_{11}(\nu)=\left(\begin{array}{cccccc}
\beta_1&0&\cdots&0\\
-i\sigma_1\nu e^{\frac{-i\nu x_0}{\sigma_1}}&i\sigma_1\nu e^{\frac{-i\nu x_1}{\sigma_1}}  &0\\
0&-i\sigma_2\nu e^{\frac{-i\nu x_1}{\sigma_2}}&i\sigma_2\nu e^{\frac{-i\nu x_2}{\sigma_2}} \\
\vdots&\hspace{.8in}\ddots&\hspace{1.2in}\ddots&\\
0&&-i\sigma_{n+1}\nu e^{\frac{-i\nu x_n}{\sigma_{n+1}}}&i\sigma_{n+1}\nu e^{\frac{-i\nu x_{n+1}}{\sigma_{n+1}}}
\end{array}\right),
$$

$$\mathcal{A}^{(p)}_{12}(\nu)=\left(\begin{array}{ccccccc}
\beta_2&0&\cdots&&0\\
-\sigma_1^2 e^{\frac{-i\nu x_0}{\sigma_1}}&\sigma_2^2  e^{\frac{-i\nu x_1}{\sigma_1}}  &0&\cdots&0\\
0&-\sigma_2^2 e^{\frac{-i\nu x_1}{\sigma_2}}&\sigma_3^2 e^{\frac{-i\nu x_2}{\sigma_2}}&0 \\
\vdots&\hspace{.8in}\ddots&\hspace{1.2in}\ddots&\hspace{1.2in}\ddots\\
0&&-\sigma_{n}^2 e^{\frac{-i\nu x_{n-1}}{\sigma_{n}}}&\sigma_{n+1}^2 e^{\frac{-i\nu x_{n}}{\sigma_{n}}}&0\\
0&&&-\sigma_{n+1}^2 e^{\frac{-i\nu x_n}{\sigma_{n+1}}}&\sigma_{n+1}^2 e^{\frac{-i\nu x_{n+1}}{\sigma_{n+1}}}
\end{array}\right),
$$
$$\mathcal{A}^{(p)}_{21}(\nu)=\left(\begin{array}{cccccc}
i\sigma_1\nu e^{\frac{i\nu x_0}{\sigma_1}}&-i\sigma_1\nu e^{\frac{i\nu x_1}{\sigma_1}}  &0\\
0&i\sigma_2\nu e^{\frac{i\nu x_1}{\sigma_2}}&-i\sigma_2\nu e^{\frac{i\nu x_2}{\sigma_2}} \\
\vdots&\hspace{.8in}\ddots&\hspace{1.2in}\ddots&\hspace{1.2in}\ddots\\
0&&i\sigma_{n+1}\nu e^{\frac{i\nu x_n}{\sigma_{n+1}}}&-i\sigma_{n+1}\nu e^{\frac{i\nu x_{n+1}}{\sigma_{n+1}}}\\
0&\cdots&0&\beta_3
\end{array}\right),
$$

$$\mathcal{A}^{(p)}_{22}(\nu)=\left(\begin{array}{ccccccc}
-\sigma_1^2 e^{\frac{i\nu x_0}{\sigma_1}}&\sigma_2^2  e^{\frac{i\nu x_1}{\sigma_1}}&  0&\cdots&0\\
0&-\sigma_2^2 e^{\frac{i\nu x_1}{\sigma_2}}&\sigma_3^2 e^{\frac{i\nu x_2}{\sigma_2}}&0 \\
\vdots&\hspace{.8in}\ddots&\hspace{1.2in}\ddots&\\
0&&-\sigma_{n}^2 e^{\frac{i\nu x_{n-1}}{\sigma_{n}}}&\sigma_{n+1}^2 e^{\frac{i\nu x_{n}}{\sigma_{n}}}&0\\
0&&&-\sigma_{n+1}^2 e^{\frac{i\nu x_n}{\sigma_{n+1}}}&\sigma_{n+1}^2 e^{\frac{i\nu x_{n+1}}{\sigma_{n+1}}}\\
0&\cdots&&0&\beta_4
\end{array}\right),
$$
}
and 
\begin{equation}\label{p_Amatrix}
\mathcal{A}^{(p)}(\nu)=\left(\begin{array}{c:c}
\mathcal{A}^{(p)}_{11}(\nu)&\mathcal{A}^{(p)}_{12}(\nu)\\
\hdashline
\mathcal{A}^{(p)}_{21}(\nu)&\mathcal{A}^{(p)}_{22}(\nu)
\end{array}
\right).
\end{equation}
\end{subequations}
and $Y(\nu,T)$, $\mathcal{Y}(\nu,T)$ are as in~\eqref{Y} and~\eqref{badY}.

The matrix $\mathcal{A}^{(p)}(\nu)$ is made up of four $(n+2)\times(n+2)$ blocks as indicated by the dashed lines.  $\mathcal{A}_{11}$ and $\mathcal{A}_{12}$ have nonzero entries only on the main and $-1$ diagonals while $\mathcal{A}_{21}$ and $\mathcal{A}_{22}$ have nonzero entries on the main and $+1$ diagonals.  The boundary conditions are incorporated in the first and last rows of $\mathcal{A}^{(p)}(\nu)$.

As before, every term in the linear equation $\mathcal{A}^{(p)}(\nu)X^{(p)}(\nu^2,T)=Y(\nu,T)$ is known.  By substituting the solutions of this equation into~\eqref{p_badsoln2}, we have solved the heat equation on the finite interval with $n$ interfaces with perfect interface conditions in terms of only known functions.  The contribution to the solution from the linear equation $\mathcal{A}^{(p)}(\nu)X^{(p)}(\nu^2,T)=\mathcal{Y}(\nu,T)$ is 0 when substituted into~\eqref{p_badsoln2} just as in Section~\ref{sec:imperfect}.

\section{Numerics}

The Fokas Method solutions presented here can be numerically implemented in a simple way using, for instance, {\sc MatLab}.  The author has done this and her code can be downloaded from \url{https://github.com/nsheils/UTM_Heat}.  A few points regarding the implementation are of note.  First, following~\cite{FlyerFokas} we parameterize $D_R^\pm$ as $\pm i\sin(\pi/8-i\theta)$ respectively.  This has the advantage of exponential decay of the integrands for both $x$ and $t$ and the points are spaced closer together near the origin.  Numerically, we only need to integrate for approximately $-10\leq \theta\leq 10$ since for $\theta$ outside this range the integrand is $0$ to machine precision.  Second, in our {\sc MatLab} implementation, we chose to scale $Y(\nu,T)$ by multiplying it by $e^{-\nu^2t}$.  This greatly improves the numerical accuracy.  Third, for large enough values of $\nu$, entries of $\mathcal{A}(\nu)$ and $\mathcal{A}^{p}(\nu)$ are too large to be represented in {\sc MatLab}.  In this case, rather than solving the linear system, we use a shape-preserving piecewise cubic Hermite polynomial to interpolate from the values for $g_\ell^{(j)}$ and $h_\ell^{(j)}$ we already have to those we still need to compute.  The code written by the author could certainly be further optimized to reduce the time it takes to run.  However, since our code confirms the accuracy of the very fast ``semi-analytical" method due to~\cite{CarrTurner}, we propose our method as a benchmark rather than a replacement.  This is especially true in the cases where their analytical method does not work ($n>9$).

In the first example we compare our solution with the exact solution.  In the second, fourth, and fifth examples we compare our solution to the one found using the code written by Carr and Turner~\cite{CarrTurner}.  A notable difference between our method and that of~\cite{CarrTurner} is in the case of time-dependent boundary conditions as in Example~\ref{exC}.  By making a transformation which makes the problem forced and alters the initial condition one can remove any time-dependence in the boundary conditions.  However, it is not clear how one would make this change in the code presented in~\cite{CarrTurner}.  In our code, the ability to use time-dependent boundary conditions is built in.  

In this manuscript we present six examples.  All times are for a 2014 MacBook Pro with a 2.8 GHz core.
\begin{enumerate}[{Example} A)]
\item\label{exA} The first example is run as a test since we can solve the problem exactly using a Fourier series solution.  We choose $n=2$, $\sigma_j=1$ for $j=1,2,3$ and the $x_j$ evenly spaced between $0$ and $1$ with perfect thermal contact.  Initially $u(x,0)=x^3$ and the boundary conditions are $u(0,t)=0$ and $u(1,t)=1.$  This example took 45.036422 seconds to run.  The relative error
\begin{equation}\label{error}
E=\frac{\max_{1\leq j\leq N}\big|u(x_j,t)-U(x_j,t)\big|}{\max_{1\leq j\leq N}\big|u(x_j,t)\big|},
\end{equation}
where $u(x_j,t)$ is the solution found using the Fokas Method method evaluated at grid points $(x_j,t)$ and $U(x_j,t)$ is the exact solution found using a Fourier series.  Our method is faulty near the end points ($x=x_0, x_{n+1}$) whenever the boundary conditions are nonhomogenous.  Thus, in the computation of the error we will omit the grid points $x_0$ and $x_{n+1}$.  However, the solution in Figure~\ref{fig:exA}  does show the Fokas Method solution evaluated at these grid points.  The results are summarized in Table~\ref{t:errorA}.

\begin{table}[h]
\caption{Relative Error in Example A}
\centering
\begin{tabular}{|l|c|c|c|}
\hline\hline
& $t=.01$ & $t=.1$&$t=1$\\ \hline
error (Fokas Method)&$3.76\times10^{-8}$ &$3.77\times10^{-8}$ &$3.75\times10^{-8}$\\
error (analytical)&$3.85\times10^{-9}$ &$3.81\times10^{-10} $&$5.16\times10^{-14}$\\
error (semi-analytical)&$3.82\times10^{-3}$ &$3.15\times10^{-3}$ &$8.94\times10^{-7}$\\
\hline
\end{tabular}
\label{t:errorA}
\end{table}
Although our method takes longer to evaluate it is much more accurate than the semi-analytical method proposed in~\cite{CarrTurner} and is similar in accuracy to the ``analytical" approach they propose.  Further, as shown in the following examples, our method works when their analytical method fails (\emph{i.e.} in the case of large $n$).

In Figure~\ref{fig:exA} the true solution is plotted as a solid line in black and the computed solution is plotted as a dashed line in red.  

\begin{figure}[htbp]
\begin{center}
\def\svgwidth{.6\textwidth}
   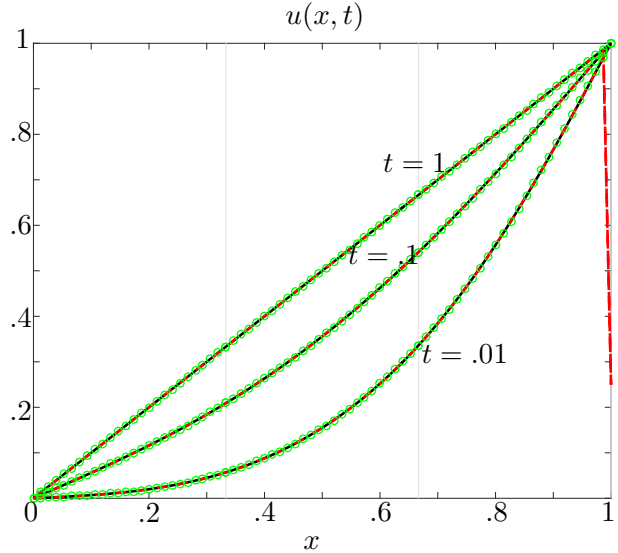 
   \caption{The true solution to the heat equation with $\sigma=1$, $u(x,0)=x^3$, $u(0,t)=0$, and $u(1,t)=1$ as described in Example~\ref{exA} is plotted as a solid line in black and the computed Fokas Method solution is plotted as a dashed line in red  \label{fig:exA}}
  \end{center}
\end{figure}

In the case when $u(x_0,t)=0$ and $u(x_{n+1},t)=0$ such as Example~\ref{exA} with $u(1,t)=0$ the relative error is computed on all grid points (including $x_0=0$ and $x_{n+1}=1$) and the error (as summarized in Table~\ref{t:errorA2}) is the same for the analytical method of~\cite{CarrTurner} and the Fokas Method.
\begin{table}[h]
\caption{Relative Error in Example A with $u(1,t)=0$.}
\centering
\begin{tabular}{|l|c|c|c|}
\hline\hline
& $t=.001$ & $t=.01$&$t=.1$\\ \hline
error (Fokas Method)&$7.06\times10^{-3}$ &$1.27\times10^{-3} $&$5.37\times10^{-4}$\\
error (analytical)&$7.06\times10^{-3} $&$1.27\times10^{-3} $&$5.37\times10^{-4}$\\
error (semi-analytical)&$9.70\times10^{-4}$ &$1.68\times10^{-3} $&$1.58\times10^{-3}$\\
\hline
\end{tabular}
\label{t:errorA2}
\end{table}

\item\label{exB} In this example we take $n=9$, $\sigma_j=1$ for $j=1,3,5,7, 9$, and $\sigma_j=\sqrt{.1}$ for $j=2,4,6,8$.  We let $u(x,0)=0$, $u(0,t)=1$ and $u(1,t)=0.$  Again, the $x_j$ are evenly spaced and we assume perfect thermal contact.  This example is the same as Example C in~\cite{CarrTurner}.  In Figure~\ref{fig:exB} our code took 75.501220 seconds to evaluate (red dashed line) while the semi-analytical method due to~\cite{CarrTurner} took 1.361871 seconds (solid blue line), and their analytical method took 2.116983 seconds (green circles).

\begin{figure}[htbp]
\begin{center}
\def\svgwidth{.6\textwidth}
   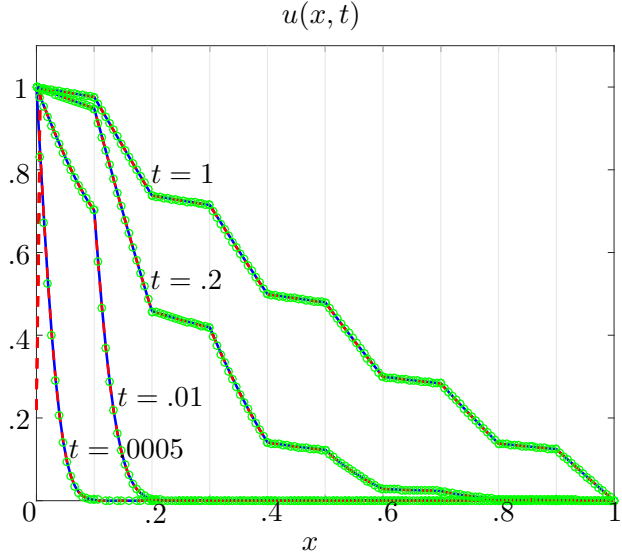 
   \caption{The solution to the heat equation with alternating diffusivities $1$ and $\sqrt{.1}$, $u(x,0)=0$, $u(0,t)=1$, and $u(1,t)=0$ with perfect thermal contact via the semi-analytical method of~\cite{CarrTurner} is plotted as a solid line in blue, their analytical method in green circles, and the Fokas Method solution is plotted as a dashed line in red as described in Example~\ref{exB}.  \label{fig:exB}}
  \end{center}
\end{figure}

\item\label{exC} In this example we take $n=3$, $\sigma_1=\sqrt{.2}$, $\sigma_2=\sqrt{.01}$, $\sigma_3=\sqrt{.1}$, and $\sigma_4=1$.  We let $u(x,0)=1$, $u(0,t)=\cos(t)$ and $u(1,t)+u_x(1,t)=0.$  Again, the $x_j$ are evenly spaced and we assume perfect thermal contact.  In Figure~\ref{fig:exC} our code took 53.407197 seconds to evaluate (red solid line).  The code provided in~\cite{CarrTurner} does not readily adapt to time-dependent boundary conditions.

\begin{figure}[htbp]
\begin{center}
\def\svgwidth{.6\textwidth}
   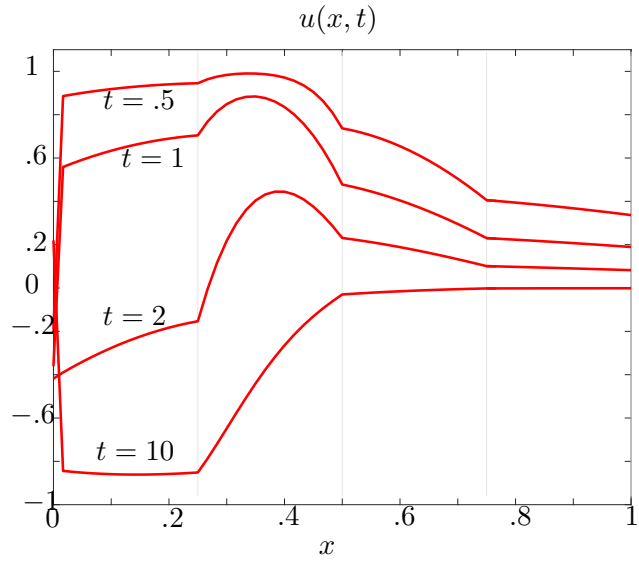 
      \caption{The solution to the heat equation with diffusivities $\sigma_1=\sqrt{.2}$, $\sigma_2=\sqrt{.01}$, $\sigma_3=\sqrt{.1}$, and $\sigma_4=1$, initial condition $u(x,0)=1$, boundary conditions $u(0,t)=\cos(t)$, and $u(1,t)+u_x(1,t)=0$, and perfect thermal contact is plotted as a solid line in red as described in Example~\ref{exC}.  \label{fig:exC}}  \end{center}
\end{figure}

\item\label{exD} This example is the same as Example B except that $u_x(1,t)=0$ and we assume imperfect thermal contact with $H_j=1/2$ for $j=1,\cdots, n$.  This example is the same as Example D in~\cite{CarrTurner}.  In Figure~\ref{fig:exC} our code took 53.092184 seconds to evaluate (red dashed line) while the semi-analytical method due to~\cite{CarrTurner}  took 0.849591 seconds (blue solid line).
\begin{figure}[htbp]
\begin{center}
\def\svgwidth{.6\textwidth}
   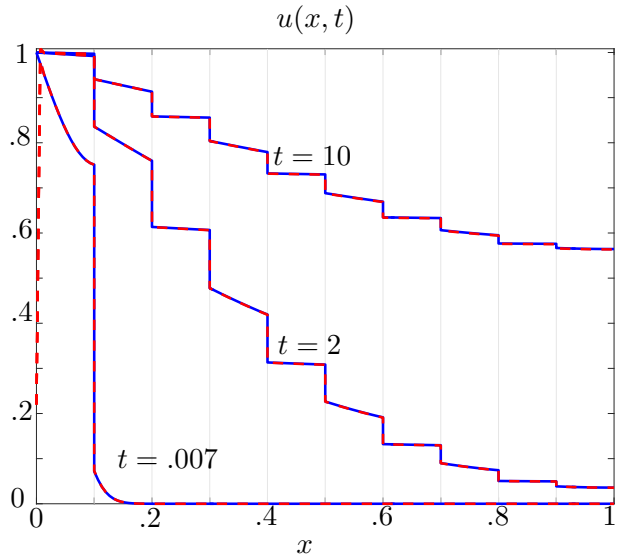 
      \caption{The solution to the heat equation with alternating diffusivities $1$ and $\sqrt{.1}$, $u(x,0)=0$, $u(0,t)=1$, and $u_x(1,t)=0$ with imperfect thermal contact and $H_j=1/2$ for $j=1,\cdots,n$ via the semi-analytical method of~\cite{CarrTurner} is plotted as a solid line in blue and the Fokas Method solution is plotted as a dashed line in red as described in Example~\ref{exD}.  \label{fig:exD}}  \end{center}
\end{figure}

\item\label{exE} For this example we take $n=199$, $\sigma_j=\sqrt{1.1+\sin(j)}$ for $j=1,\cdots, n$.  We let $u(x,0)=1$, $u(0,t)=1/2$ and $u(1,t)=0.$  Again, the $x_j$ are evenly spaced and we assume perfect thermal contact.  This example is the same as the macroscopic modeling example in~\cite{CarrTurner}.  In Figure~\ref{fig:exE} our code took 2067.976494 seconds to evaluate (red dashed line) while the semi-analytical method due to~\cite{CarrTurner} took 8.078391 seconds (blue solid line).
\begin{figure}[htbp]
\begin{center}
\def\svgwidth{.6\textwidth}
   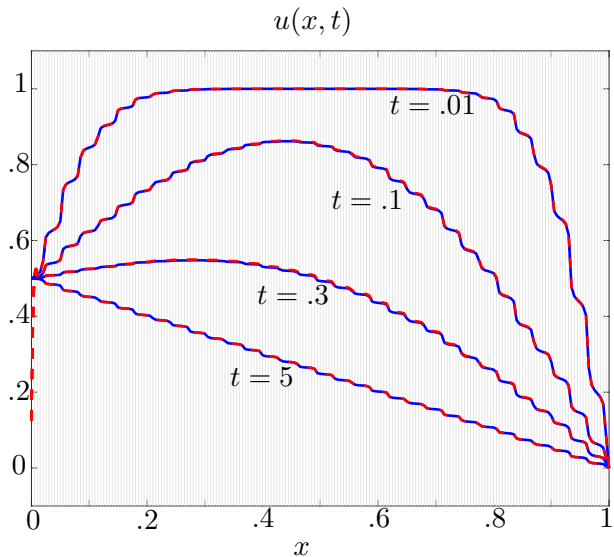 
      \caption{The solution to the heat equation with diffusivities $\sigma_j=\sqrt{1.1+\sin(j)}$ for $j=1,\cdots, n$, $u(x,0)=1$, $u(0,t)=1/2$, and $u(1,t)=0$ with perfect thermal contact via the semi-analytical method of~\cite{CarrTurner} is plotted as a solid line in blue and the Fokas Method solution is plotted as a dashed line in red as described in Example~\ref{exE}.  \label{fig:exE}}  \end{center}
\end{figure}

\item\label{exF}  In our final example we take $n=199$, $\sigma_j=\sqrt{1.1+\sin(j)}$ for $j=1,\cdots, n$.  We let $u(x,0)=x$, $u_x(0,t)=0$ and $u(1,t)=0$  The $x_j$ are evenly spaced and we assume imperfect thermal contact with $H_j=1/2$ for $j=1,\cdots,n$.  In Figure~\ref{fig:exF} our code took 1064.144519 seconds to evaluate (red dashed line) while the semi-analytical method due to~\cite{CarrTurner} took 5.075486 seconds (solid blue line).
\begin{figure}[htbp]
\begin{center}
\def\svgwidth{.6\textwidth}
   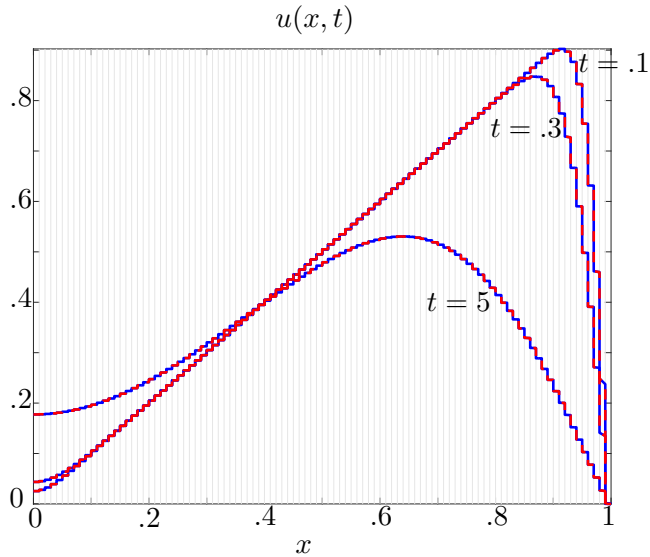 
      \caption{The solution to the heat equation with diffusivities $\sigma_j=\sqrt{1.1+\sin(j)}$ for $j=1,\cdots, n$, $u(x,0)=x$, $u_x(0,t)=0$, $u(1,t)=0$, $H_j=1/2$ for $j=1,\cdots,n$ with imperfect thermal contact via the semi-analytical method of~\cite{CarrTurner} is plotted as a solid line in blue and the Fokas Method solution is plotted as a dashed line in red as described in Example~\ref{exF}.  \label{fig:exF}}  \end{center}
\end{figure}

\end{enumerate}

\section{Conclusion}
In this manuscript the Fokas Method is to provide explicit solution formulae for the heat transport interface problem with perfect and imperfect interface conditions with arbitrary boundary conditions and a generic number of interfaces.  The generality of the method would also easily allow for other possible interface conditions.  Further, numerical implementation of the solutions are provided.  Although the proposed code is slower than other possibilities, and could certainly be further optimized, it provides a good option for benchmarking other schemes which rely on less explicit analytical solutions.

\section{Acknowledgements}
I would like to acknowledge Bernard Deconinck for bringing this problem to my attention and for useful conversations.  I would also like to thank Elliot Carr and Ian Turner for making their well-documented code available.  Olga Trichtchenko and Christopher Swierczewski also provided useful feedback.

\newpage
\bibliographystyle{abbrv}
\bibliography{FullBib}

\end{document}

%% file: GR_domain_nf.pdf_tex
\begingroup%
  \makeatletter%
  \providecommand\color[2][]{%
    \errmessage{(Inkscape) Color is used for the text in Inkscape, but the package 'color.sty' is not loaded}%
    \renewcommand\color[2][]{}%
  }%
  \providecommand\transparent[1]{%
    \errmessage{(Inkscape) Transparency is used (non-zero) for the text in Inkscape, but the package 'transparent.sty' is not loaded}%
    \renewcommand\transparent[1]{}%
  }%
  \providecommand\rotatebox[2]{#2}%
  \ifx\svgwidth\undefined%
    \setlength{\unitlength}{1067.857664bp}%
    \ifx\svgscale\undefined%
      \relax%
    \else%
      \setlength{\unitlength}{\unitlength * \real{\svgscale}}%
    \fi%
  \else%
    \setlength{\unitlength}{\svgwidth}%
  \fi%
  \global\let\svgwidth\undefined%
  \global\let\svgscale\undefined%
  \makeatother%
  \begin{picture}(1,0.29412871)%
    \put(0,0){\includegraphics[width=\unitlength]{GR_domain_nf.pdf}}%
    \put(0.97269889,0.04427738){\color[rgb]{0,0,0}\makebox(0,0)[lt]{\begin{minipage}{0.03578399\unitlength}\raggedright $x$\end{minipage}}}%
    \put(0.16236052,0.02188992){\color[rgb]{0,0,0}\makebox(0,0)[lt]{\begin{minipage}{0.11490773\unitlength}\raggedright $x_1$\end{minipage}}}%
    \put(0.91326297,0.02210469){\color[rgb]{0,0,0}\makebox(0,0)[lt]{\begin{minipage}{0.10834471\unitlength}\raggedright $x_{n+1}$\end{minipage}}}%
    \put(0.01072965,0.02188868){\color[rgb]{0,0,0}\makebox(0,0)[lt]{\begin{minipage}{0.09284909\unitlength}\raggedright $x_0=0$\end{minipage}}}%
    \put(0.0112459,0.29709913){\color[rgb]{0,0,0}\makebox(0,0)[lt]{\begin{minipage}{0.0348423\unitlength}\raggedright $t$\end{minipage}}}%
    \put(0.02030978,0.24327065){\color[rgb]{0,0,0}\makebox(0,0)[lt]{\begin{minipage}{0.10980415\unitlength}\raggedright $T$\end{minipage}}}%
    \put(0.76139965,0.02188091){\color[rgb]{0,0,0}\makebox(0,0)[lt]{\begin{minipage}{0.08583226\unitlength}\raggedright $x_n$\end{minipage}}}%
    \put(0.61088224,0.02209554){\color[rgb]{0,0,0}\makebox(0,0)[lt]{\begin{minipage}{0.1174547\unitlength}\raggedright $x_{n-1}$\end{minipage}}}%
    \put(0.3073034,0.02188274){\color[rgb]{0,0,0}\makebox(0,0)[lt]{\begin{minipage}{0.10480662\unitlength}\raggedright $x_2$\end{minipage}}}%
    \put(0.44096983,0.01813971){\color[rgb]{0,0,0}\makebox(0,0)[lt]{\begin{minipage}{0.10632557\unitlength}\raggedright $\cdots$\end{minipage}}}%
  \end{picture}%
\endgroup%

%% file: heat_Dpm.pdf_tex
\begingroup%
  \makeatletter%
  \providecommand\color[2][]{%
    \errmessage{(Inkscape) Color is used for the text in Inkscape, but the package 'color.sty' is not loaded}%
    \renewcommand\color[2][]{}%
  }%
  \providecommand\transparent[1]{%
    \errmessage{(Inkscape) Transparency is used (non-zero) for the text in Inkscape, but the package 'transparent.sty' is not loaded}%
    \renewcommand\transparent[1]{}%
  }%
  \providecommand\rotatebox[2]{#2}%
  \ifx\svgwidth\undefined%
    \setlength{\unitlength}{535.09785813bp}%
    \ifx\svgscale\undefined%
      \relax%
    \else%
      \setlength{\unitlength}{\unitlength * \real{\svgscale}}%
    \fi%
  \else%
    \setlength{\unitlength}{\svgwidth}%
  \fi%
  \global\let\svgwidth\undefined%
  \global\let\svgscale\undefined%
  \makeatother%
  \begin{picture}(1,0.94782192)%
    \put(0,0){\includegraphics[width=\unitlength,page=1]{heat_Dpm.pdf}}%
    \put(0.4321468,0.96199277){\color[rgb]{0,0,0}\makebox(0,0)[lt]{\begin{minipage}{0.33788212\unitlength}\raggedright $\Im(\nu)$\end{minipage}}}%
    \put(0.92252439,0.4805855){\color[rgb]{0,0,0}\makebox(0,0)[lt]{\begin{minipage}{0.36443538\unitlength}\raggedright $\Re(\nu)$\end{minipage}}}%
    \put(0.41214063,0.74696077){\color[rgb]{0,0,0}\makebox(0,0)[lt]{\begin{minipage}{0.36166528\unitlength}\raggedright $D^+$\end{minipage}}}%
    \put(0,0){\includegraphics[width=\unitlength,page=2]{heat_Dpm.pdf}}%
    \put(0.46130035,0.22966566){\color[rgb]{0,0,0}\makebox(0,0)[lt]{\begin{minipage}{0.33863622\unitlength}\raggedright $D^-$\end{minipage}}}%
    \put(0,0){\includegraphics[width=\unitlength,page=3]{heat_Dpm.pdf}}%
  \end{picture}%
\endgroup%

%% file: heat_DRpm.pdf_tex
\begingroup%
  \makeatletter%
  \providecommand\color[2][]{%
    \errmessage{(Inkscape) Color is used for the text in Inkscape, but the package 'color.sty' is not loaded}%
    \renewcommand\color[2][]{}%
  }%
  \providecommand\transparent[1]{%
    \errmessage{(Inkscape) Transparency is used (non-zero) for the text in Inkscape, but the package 'transparent.sty' is not loaded}%
    \renewcommand\transparent[1]{}%
  }%
  \providecommand\rotatebox[2]{#2}%
  \ifx\svgwidth\undefined%
    \setlength{\unitlength}{535.09785813bp}%
    \ifx\svgscale\undefined%
      \relax%
    \else%
      \setlength{\unitlength}{\unitlength * \real{\svgscale}}%
    \fi%
  \else%
    \setlength{\unitlength}{\svgwidth}%
  \fi%
  \global\let\svgwidth\undefined%
  \global\let\svgscale\undefined%
  \makeatother%
  \begin{picture}(1,0.94782192)%
    \put(0,0){\includegraphics[width=\unitlength,page=1]{heat_DRpm.pdf}}%
    \put(0.4321468,0.96199277){\color[rgb]{0,0,0}\makebox(0,0)[lt]{\begin{minipage}{0.33788212\unitlength}\raggedright $\Im(\nu)$\end{minipage}}}%
    \put(0.92252439,0.4805855){\color[rgb]{0,0,0}\makebox(0,0)[lt]{\begin{minipage}{0.36443538\unitlength}\raggedright $\Re(\nu)$\end{minipage}}}%
    \put(0.43008127,0.76789152){\color[rgb]{0,0,0}\makebox(0,0)[lt]{\begin{minipage}{0.36166528\unitlength}\raggedright $D_R^+$\end{minipage}}}%
    \put(0,0){\includegraphics[width=\unitlength,page=2]{heat_DRpm.pdf}}%
    \put(0.47027067,0.22069534){\color[rgb]{0,0,0}\makebox(0,0)[lt]{\begin{minipage}{0.33863622\unitlength}\raggedright $D_R^-$\end{minipage}}}%
    \put(0,0){\includegraphics[width=\unitlength,page=3]{heat_DRpm.pdf}}%
    \put(0.38621968,0.53664984){\color[rgb]{0,0,0}\makebox(0,0)[lt]{\begin{minipage}{0.17626161\unitlength}\raggedright $R$\end{minipage}}}%
    \put(0,0){\includegraphics[width=\unitlength,page=4]{heat_DRpm.pdf}}%
  \end{picture}%
\endgroup%

%% file: heat_DRpm_close.pdf_tex
\begingroup%
  \makeatletter%
  \providecommand\color[2][]{%
    \errmessage{(Inkscape) Color is used for the text in Inkscape, but the package 'color.sty' is not loaded}%
    \renewcommand\color[2][]{}%
  }%
  \providecommand\transparent[1]{%
    \errmessage{(Inkscape) Transparency is used (non-zero) for the text in Inkscape, but the package 'transparent.sty' is not loaded}%
    \renewcommand\transparent[1]{}%
  }%
  \providecommand\rotatebox[2]{#2}%
  \ifx\svgwidth\undefined%
    \setlength{\unitlength}{535.09785813bp}%
    \ifx\svgscale\undefined%
      \relax%
    \else%
      \setlength{\unitlength}{\unitlength * \real{\svgscale}}%
    \fi%
  \else%
    \setlength{\unitlength}{\svgwidth}%
  \fi%
  \global\let\svgwidth\undefined%
  \global\let\svgscale\undefined%
  \makeatother%
  \begin{picture}(1,0.94782192)%
    \put(0,0){\includegraphics[width=\unitlength,page=1]{heat_DRpm_close.pdf}}%
    \put(0.4321468,0.96199277){\color[rgb]{0,0,0}\makebox(0,0)[lt]{\begin{minipage}{0.33788212\unitlength}\raggedright $\Im(\nu)$\end{minipage}}}%
    \put(0.92252439,0.4805855){\color[rgb]{0,0,0}\makebox(0,0)[lt]{\begin{minipage}{0.36443538\unitlength}\raggedright $\Re(\nu)$\end{minipage}}}%
    \put(0,0){\includegraphics[width=\unitlength,page=2]{heat_DRpm_close.pdf}}%
    \put(0.38621968,0.53664984){\color[rgb]{0,0,0}\makebox(0,0)[lt]{\begin{minipage}{0.17626161\unitlength}\raggedright $R$\end{minipage}}}%
    \put(0,0){\includegraphics[width=\unitlength,page=3]{heat_DRpm_close.pdf}}%
    \put(0.64867481,0.65552289){\color[rgb]{0,0,0}\makebox(0,0)[lt]{\begin{minipage}{0.39575728\unitlength}\raggedright $\mathcal{L}_{D^+}$\end{minipage}}}%
    \put(0.59280924,0.77025919){\color[rgb]{0,0,0}\makebox(0,0)[lt]{\begin{minipage}{0.20951853\unitlength}\raggedright $C$\end{minipage}}}%
    \put(0.13152989,0.89196744){\color[rgb]{0,0,0}\makebox(0,0)[lt]{\begin{minipage}{0.44564259\unitlength}\raggedright $\mathcal{L}_{C}^+$\end{minipage}}}%
    \put(0,0){\includegraphics[width=\unitlength,page=4]{heat_DRpm_close.pdf}}%
    \put(0.66463209,0.28851533){\color[rgb]{0,0,0}\makebox(0,0)[lt]{\begin{minipage}{0.39575728\unitlength}\raggedright $\mathcal{L}_{D^-}$\end{minipage}}}%
    \put(0.17377219,0.13751688){\color[rgb]{0,0,0}\makebox(0,0)[lt]{\begin{minipage}{0.44564259\unitlength}\raggedright $\mathcal{L}_{C}^-$\end{minipage}}}%
    \put(0.59280924,0.19373313){\color[rgb]{0,0,0}\makebox(0,0)[lt]{\begin{minipage}{0.20951853\unitlength}\raggedright $C$\end{minipage}}}%
  \end{picture}%
\endgroup%

%% file: ExA_paper.pdf_tex
\begingroup%
  \makeatletter%
  \providecommand\color[2][]{%
    \errmessage{(Inkscape) Color is used for the text in Inkscape, but the package 'color.sty' is not loaded}%
    \renewcommand\color[2][]{}%
  }%
  \providecommand\transparent[1]{%
    \errmessage{(Inkscape) Transparency is used (non-zero) for the text in Inkscape, but the package 'transparent.sty' is not loaded}%
    \renewcommand\transparent[1]{}%
  }%
  \providecommand\rotatebox[2]{#2}%
  \ifx\svgwidth\undefined%
    \setlength{\unitlength}{576.00001843bp}%
    \ifx\svgscale\undefined%
      \relax%
    \else%
      \setlength{\unitlength}{\unitlength * \real{\svgscale}}%
    \fi%
  \else%
    \setlength{\unitlength}{\svgwidth}%
  \fi%
  \global\let\svgwidth\undefined%
  \global\let\svgscale\undefined%
  \makeatother%
  \begin{picture}(1,0.75043566)%
    \put(0,0){\includegraphics[width=\unitlength,page=1]{ExA_paper.pdf}}%
    \put(0.10621695,0.71280355){\color[rgb]{0,0,0}\makebox(0,0)[lt]{\begin{minipage}{0.07699116\unitlength}\raggedright $1$\end{minipage}}}%
    \put(0.59969863,0.54573301){\color[rgb]{0,0,0}\makebox(0,0)[lt]{\begin{minipage}{0.11814159\unitlength}\raggedright $t=1$\end{minipage}}}%
    \put(0.55250099,0.42311745){\color[rgb]{0,0,0}\makebox(0,0)[lt]{\begin{minipage}{0.16725663\unitlength}\raggedright $t=.1$\end{minipage}}}%
    \put(0.65134564,0.28966123){\color[rgb]{0,0,0}\makebox(0,0)[lt]{\begin{minipage}{0.14469023\unitlength}\raggedright $t=.01$\end{minipage}}}%
    \put(0.49252068,0.03177039){\color[rgb]{0,0,0}\makebox(0,0)[lt]{\begin{minipage}{0.28407079\unitlength}\raggedright $x$\end{minipage}}}%
    \put(0.11732811,0.07807935){\color[rgb]{0,0,0}\makebox(0,0)[lt]{\begin{minipage}{0.06106193\unitlength}\raggedright $0$\end{minipage}}}%
    \put(0.2701057,0.07575131){\color[rgb]{0,0,0}\makebox(0,0)[lt]{\begin{minipage}{0.07035399\unitlength}\raggedright $.2$\end{minipage}}}%
    \put(0.42566341,0.07808105){\color[rgb]{0,0,0}\makebox(0,0)[lt]{\begin{minipage}{0.079646\unitlength}\raggedright $.4$\end{minipage}}}%
    \put(0.58121655,0.07807935){\color[rgb]{0,0,0}\makebox(0,0)[lt]{\begin{minipage}{0.0637168\unitlength}\raggedright $.6$\end{minipage}}}%
    \put(0.73677388,0.07808113){\color[rgb]{0,0,0}\makebox(0,0)[lt]{\begin{minipage}{0.08628319\unitlength}\raggedright $.8$\end{minipage}}}%
    \put(0.89232968,0.0780802){\color[rgb]{0,0,0}\makebox(0,0)[lt]{\begin{minipage}{0.07168138\unitlength}\raggedright $1$\end{minipage}}}%
    \put(0.09788343,0.21280284){\color[rgb]{0,0,0}\makebox(0,0)[lt]{\begin{minipage}{0.08893805\unitlength}\raggedright $.2$\end{minipage}}}%
    \put(0.09788391,0.33780428){\color[rgb]{0,0,0}\makebox(0,0)[lt]{\begin{minipage}{0.07168141\unitlength}\raggedright $.4$\end{minipage}}}%
    \put(0.09788549,0.46280347){\color[rgb]{0,0,0}\makebox(0,0)[lt]{\begin{minipage}{0.08362832\unitlength}\raggedright $.6$\end{minipage}}}%
    \put(0.09788436,0.58780207){\color[rgb]{0,0,0}\makebox(0,0)[lt]{\begin{minipage}{0.09159293\unitlength}\raggedright $.8$\end{minipage}}}%
    \put(0.46862684,0.74858452){\color[rgb]{0,0,0}\makebox(0,0)[lt]{\begin{minipage}{0.37433627\unitlength}\raggedright $u(x,t)$\end{minipage}}}%
  \end{picture}%
\endgroup%

%% file: ExB_paper.pdf_tex
\begingroup%
  \makeatletter%
  \providecommand\color[2][]{%
    \errmessage{(Inkscape) Color is used for the text in Inkscape, but the package 'color.sty' is not loaded}%
    \renewcommand\color[2][]{}%
  }%
  \providecommand\transparent[1]{%
    \errmessage{(Inkscape) Transparency is used (non-zero) for the text in Inkscape, but the package 'transparent.sty' is not loaded}%
    \renewcommand\transparent[1]{}%
  }%
  \providecommand\rotatebox[2]{#2}%
  \ifx\svgwidth\undefined%
    \setlength{\unitlength}{576.00001843bp}%
    \ifx\svgscale\undefined%
      \relax%
    \else%
      \setlength{\unitlength}{\unitlength * \real{\svgscale}}%
    \fi%
  \else%
    \setlength{\unitlength}{\svgwidth}%
  \fi%
  \global\let\svgwidth\undefined%
  \global\let\svgscale\undefined%
  \makeatother%
  \begin{picture}(1,0.75)%
    \put(0,0){\includegraphics[width=\unitlength,page=1]{ExB_paper.pdf}}%
    \put(0.09997945,0.64774706){\color[rgb]{0,0,0}\makebox(0,0)[lt]{\begin{minipage}{0.07699116\unitlength}\raggedright $1$\end{minipage}}}%
    \put(0.2840462,0.53382303){\color[rgb]{0,0,0}\makebox(0,0)[lt]{\begin{minipage}{0.11814159\unitlength}\raggedright $t=1$\end{minipage}}}%
    \put(0.28343164,0.39117305){\color[rgb]{0,0,0}\makebox(0,0)[lt]{\begin{minipage}{0.16725663\unitlength}\raggedright $t=.2$\end{minipage}}}%
    \put(0.23903636,0.23431467){\color[rgb]{0,0,0}\makebox(0,0)[lt]{\begin{minipage}{0.14469023\unitlength}\raggedright $t=.01$\end{minipage}}}%
    \put(0.48628318,0.03338056){\color[rgb]{0,0,0}\makebox(0,0)[lt]{\begin{minipage}{0.28407079\unitlength}\raggedright $x$\end{minipage}}}%
    \put(0.11109061,0.07968953){\color[rgb]{0,0,0}\makebox(0,0)[lt]{\begin{minipage}{0.06106193\unitlength}\raggedright $0$\end{minipage}}}%
    \put(0.27497931,0.07736148){\color[rgb]{0,0,0}\makebox(0,0)[lt]{\begin{minipage}{0.07035399\unitlength}\raggedright $.2$\end{minipage}}}%
    \put(0.43053703,0.07969122){\color[rgb]{0,0,0}\makebox(0,0)[lt]{\begin{minipage}{0.079646\unitlength}\raggedright $.4$\end{minipage}}}%
    \put(0.58609016,0.07968953){\color[rgb]{0,0,0}\makebox(0,0)[lt]{\begin{minipage}{0.0637168\unitlength}\raggedright $.6$\end{minipage}}}%
    \put(0.74164749,0.07969131){\color[rgb]{0,0,0}\makebox(0,0)[lt]{\begin{minipage}{0.08628319\unitlength}\raggedright $.8$\end{minipage}}}%
    \put(0.8972033,0.07969038){\color[rgb]{0,0,0}\makebox(0,0)[lt]{\begin{minipage}{0.07168138\unitlength}\raggedright $1$\end{minipage}}}%
    \put(0.09164593,0.20330191){\color[rgb]{0,0,0}\makebox(0,0)[lt]{\begin{minipage}{0.08893805\unitlength}\raggedright $.2$\end{minipage}}}%
    \put(0.09164641,0.31163668){\color[rgb]{0,0,0}\makebox(0,0)[lt]{\begin{minipage}{0.07168141\unitlength}\raggedright $.4$\end{minipage}}}%
    \put(0.09164799,0.42552476){\color[rgb]{0,0,0}\makebox(0,0)[lt]{\begin{minipage}{0.08362832\unitlength}\raggedright $.6$\end{minipage}}}%
    \put(0.09164686,0.53385669){\color[rgb]{0,0,0}\makebox(0,0)[lt]{\begin{minipage}{0.09159293\unitlength}\raggedright $.8$\end{minipage}}}%
    \put(0.45961156,0.75297247){\color[rgb]{0,0,0}\makebox(0,0)[lt]{\begin{minipage}{0.37433627\unitlength}\raggedright $u(x,t)$\end{minipage}}}%
    \put(0.17384461,0.16467357){\color[rgb]{0,0,0}\makebox(0,0)[lt]{\begin{minipage}{0.16327429\unitlength}\raggedright $t=.0005$\end{minipage}}}%
  \end{picture}%
\endgroup%

%% file: ExC_paper.pdf_tex
\begingroup%
  \makeatletter%
  \providecommand\color[2][]{%
    \errmessage{(Inkscape) Color is used for the text in Inkscape, but the package 'color.sty' is not loaded}%
    \renewcommand\color[2][]{}%
  }%
  \providecommand\transparent[1]{%
    \errmessage{(Inkscape) Transparency is used (non-zero) for the text in Inkscape, but the package 'transparent.sty' is not loaded}%
    \renewcommand\transparent[1]{}%
  }%
  \providecommand\rotatebox[2]{#2}%
  \ifx\svgwidth\undefined%
    \setlength{\unitlength}{576.00001843bp}%
    \ifx\svgscale\undefined%
      \relax%
    \else%
      \setlength{\unitlength}{\unitlength * \real{\svgscale}}%
    \fi%
  \else%
    \setlength{\unitlength}{\svgwidth}%
  \fi%
  \global\let\svgwidth\undefined%
  \global\let\svgscale\undefined%
  \makeatother%
  \begin{picture}(1,0.75)%
    \put(0,0){\includegraphics[width=\unitlength,page=1]{ExC_paper.pdf}}%
    \put(0.22165681,0.56157831){\color[rgb]{0,0,0}\makebox(0,0)[lt]{\begin{minipage}{0.11814159\unitlength}\raggedright $t=1$\end{minipage}}}%
    \put(0.19714845,0.63928232){\color[rgb]{0,0,0}\makebox(0,0)[lt]{\begin{minipage}{0.16725663\unitlength}\raggedright $t=.5$\end{minipage}}}%
    \put(0.19655849,0.34859896){\color[rgb]{0,0,0}\makebox(0,0)[lt]{\begin{minipage}{0.14469023\unitlength}\raggedright $t=2$\end{minipage}}}%
    \put(0.48628318,0.03350556){\color[rgb]{0,0,0}\makebox(0,0)[lt]{\begin{minipage}{0.28407079\unitlength}\raggedright $x$\end{minipage}}}%
    \put(0.11942394,0.07745342){\color[rgb]{0,0,0}\makebox(0,0)[lt]{\begin{minipage}{0.06106193\unitlength}\raggedright $0$\end{minipage}}}%
    \put(0.27497931,0.07748648){\color[rgb]{0,0,0}\makebox(0,0)[lt]{\begin{minipage}{0.07035399\unitlength}\raggedright $.2$\end{minipage}}}%
    \put(0.43053703,0.07981622){\color[rgb]{0,0,0}\makebox(0,0)[lt]{\begin{minipage}{0.079646\unitlength}\raggedright $.4$\end{minipage}}}%
    \put(0.58609016,0.07981453){\color[rgb]{0,0,0}\makebox(0,0)[lt]{\begin{minipage}{0.0637168\unitlength}\raggedright $.6$\end{minipage}}}%
    \put(0.74164749,0.07981631){\color[rgb]{0,0,0}\makebox(0,0)[lt]{\begin{minipage}{0.08628319\unitlength}\raggedright $.8$\end{minipage}}}%
    \put(0.8972033,0.07981538){\color[rgb]{0,0,0}\makebox(0,0)[lt]{\begin{minipage}{0.07168138\unitlength}\raggedright $1$\end{minipage}}}%
    \put(0.08886815,0.10064913){\color[rgb]{0,0,0}\makebox(0,0)[lt]{\begin{minipage}{0.08893805\unitlength}\raggedright $-1$\end{minipage}}}%
    \put(0.07220196,0.21731724){\color[rgb]{0,0,0}\makebox(0,0)[lt]{\begin{minipage}{0.08893804\unitlength}\raggedright $-.6$\end{minipage}}}%
    \put(0.07220355,0.33676087){\color[rgb]{0,0,0}\makebox(0,0)[lt]{\begin{minipage}{0.08362832\unitlength}\raggedright $-.2$\end{minipage}}}%
    \put(0.45961156,0.75309747){\color[rgb]{0,0,0}\makebox(0,0)[lt]{\begin{minipage}{0.37433627\unitlength}\raggedright $u(x,t)$\end{minipage}}}%
    \put(0.1877335,0.16757635){\color[rgb]{0,0,0}\makebox(0,0)[lt]{\begin{minipage}{0.16327429\unitlength}\raggedright $t=10$\end{minipage}}}%
    \put(0.09164799,0.39231643){\color[rgb]{0,0,0}\makebox(0,0)[lt]{\begin{minipage}{0.08362832\unitlength}\raggedright $0$\end{minipage}}}%
    \put(0.09164593,0.44787134){\color[rgb]{0,0,0}\makebox(0,0)[lt]{\begin{minipage}{0.08893805\unitlength}\raggedright $.2$\end{minipage}}}%
    \put(0.09164641,0.56453945){\color[rgb]{0,0,0}\makebox(0,0)[lt]{\begin{minipage}{0.07168141\unitlength}\raggedright $.6$\end{minipage}}}%
    \put(0.09164799,0.68398309){\color[rgb]{0,0,0}\makebox(0,0)[lt]{\begin{minipage}{0.08362832\unitlength}\raggedright $1$\end{minipage}}}%
  \end{picture}%
\endgroup%

%% file: ExD_paper.pdf_tex
\begingroup%
  \makeatletter%
  \providecommand\color[2][]{%
    \errmessage{(Inkscape) Color is used for the text in Inkscape, but the package 'color.sty' is not loaded}%
    \renewcommand\color[2][]{}%
  }%
  \providecommand\transparent[1]{%
    \errmessage{(Inkscape) Transparency is used (non-zero) for the text in Inkscape, but the package 'transparent.sty' is not loaded}%
    \renewcommand\transparent[1]{}%
  }%
  \providecommand\rotatebox[2]{#2}%
  \ifx\svgwidth\undefined%
    \setlength{\unitlength}{576.00001843bp}%
    \ifx\svgscale\undefined%
      \relax%
    \else%
      \setlength{\unitlength}{\unitlength * \real{\svgscale}}%
    \fi%
  \else%
    \setlength{\unitlength}{\svgwidth}%
  \fi%
  \global\let\svgwidth\undefined%
  \global\let\svgscale\undefined%
  \makeatother%
  \begin{picture}(1,0.75)%
    \put(0,0){\includegraphics[width=\unitlength,page=1]{ExD_paper.pdf}}%
    \put(0.45425301,0.30855509){\color[rgb]{0,0,0}\makebox(0,0)[lt]{\begin{minipage}{0.11814159\unitlength}\raggedright $t=2$\end{minipage}}}%
    \put(0.44788625,0.56219322){\color[rgb]{0,0,0}\makebox(0,0)[lt]{\begin{minipage}{0.16725663\unitlength}\raggedright $t=10$\end{minipage}}}%
    \put(0.47903179,0.03279306){\color[rgb]{0,0,0}\makebox(0,0)[lt]{\begin{minipage}{0.28407079\unitlength}\raggedright $x$\end{minipage}}}%
    \put(0.12050589,0.07674092){\color[rgb]{0,0,0}\makebox(0,0)[lt]{\begin{minipage}{0.06106193\unitlength}\raggedright $0$\end{minipage}}}%
    \put(0.26772792,0.07677398){\color[rgb]{0,0,0}\makebox(0,0)[lt]{\begin{minipage}{0.07035399\unitlength}\raggedright $.2$\end{minipage}}}%
    \put(0.42328564,0.07910372){\color[rgb]{0,0,0}\makebox(0,0)[lt]{\begin{minipage}{0.079646\unitlength}\raggedright $.4$\end{minipage}}}%
    \put(0.57883877,0.07910203){\color[rgb]{0,0,0}\makebox(0,0)[lt]{\begin{minipage}{0.0637168\unitlength}\raggedright $.6$\end{minipage}}}%
    \put(0.7343961,0.07910381){\color[rgb]{0,0,0}\makebox(0,0)[lt]{\begin{minipage}{0.08628319\unitlength}\raggedright $.8$\end{minipage}}}%
    \put(0.88995191,0.07910288){\color[rgb]{0,0,0}\makebox(0,0)[lt]{\begin{minipage}{0.07168138\unitlength}\raggedright $1$\end{minipage}}}%
    \put(0.09550565,0.10271441){\color[rgb]{0,0,0}\makebox(0,0)[lt]{\begin{minipage}{0.08893805\unitlength}\raggedright $0$\end{minipage}}}%
    \put(0.09550613,0.21938251){\color[rgb]{0,0,0}\makebox(0,0)[lt]{\begin{minipage}{0.08893804\unitlength}\raggedright $.2$\end{minipage}}}%
    \put(0.45236018,0.75238497){\color[rgb]{0,0,0}\makebox(0,0)[lt]{\begin{minipage}{0.37433627\unitlength}\raggedright $u(x,t)$\end{minipage}}}%
    \put(0.23888919,0.15528568){\color[rgb]{0,0,0}\makebox(0,0)[lt]{\begin{minipage}{0.16327429\unitlength}\raggedright $t=.007$\end{minipage}}}%
    \put(0.09550771,0.35271504){\color[rgb]{0,0,0}\makebox(0,0)[lt]{\begin{minipage}{0.08362832\unitlength}\raggedright $.4$\end{minipage}}}%
    \put(0.09550565,0.46382551){\color[rgb]{0,0,0}\makebox(0,0)[lt]{\begin{minipage}{0.08893805\unitlength}\raggedright $.6$\end{minipage}}}%
    \put(0.09272835,0.58882695){\color[rgb]{0,0,0}\makebox(0,0)[lt]{\begin{minipage}{0.12345132\unitlength}\raggedright $.8$\end{minipage}}}%
    \put(0.10106327,0.69993725){\color[rgb]{0,0,0}\makebox(0,0)[lt]{\begin{minipage}{0.08362832\unitlength}\raggedright $1$\end{minipage}}}%
  \end{picture}%
\endgroup%

%% file: ExE_paper.pdf_tex
\begingroup%
  \makeatletter%
  \providecommand\color[2][]{%
    \errmessage{(Inkscape) Color is used for the text in Inkscape, but the package 'color.sty' is not loaded}%
    \renewcommand\color[2][]{}%
  }%
  \providecommand\transparent[1]{%
    \errmessage{(Inkscape) Transparency is used (non-zero) for the text in Inkscape, but the package 'transparent.sty' is not loaded}%
    \renewcommand\transparent[1]{}%
  }%
  \providecommand\rotatebox[2]{#2}%
  \ifx\svgwidth\undefined%
    \setlength{\unitlength}{576.00001843bp}%
    \ifx\svgscale\undefined%
      \relax%
    \else%
      \setlength{\unitlength}{\unitlength * \real{\svgscale}}%
    \fi%
  \else%
    \setlength{\unitlength}{\svgwidth}%
  \fi%
  \global\let\svgwidth\undefined%
  \global\let\svgscale\undefined%
  \makeatother%
  \begin{picture}(1,0.75)%
    \put(0,0){\includegraphics[width=\unitlength,page=1]{ExE_paper.pdf}}%
    \put(0.39702587,0.26811753){\color[rgb]{0,0,0}\makebox(0,0)[lt]{\begin{minipage}{0.11814159\unitlength}\raggedright $t=5$\end{minipage}}}%
    \put(0.61116058,0.63360424){\color[rgb]{0,0,0}\makebox(0,0)[lt]{\begin{minipage}{0.16725663\unitlength}\raggedright $t=.01$\end{minipage}}}%
    \put(0.48180957,0.03279306){\color[rgb]{0,0,0}\makebox(0,0)[lt]{\begin{minipage}{0.28407079\unitlength}\raggedright $x$\end{minipage}}}%
    \put(0.12328366,0.07674092){\color[rgb]{0,0,0}\makebox(0,0)[lt]{\begin{minipage}{0.06106193\unitlength}\raggedright $0$\end{minipage}}}%
    \put(0.2705057,0.07677398){\color[rgb]{0,0,0}\makebox(0,0)[lt]{\begin{minipage}{0.07035399\unitlength}\raggedright $.2$\end{minipage}}}%
    \put(0.42606341,0.07910372){\color[rgb]{0,0,0}\makebox(0,0)[lt]{\begin{minipage}{0.079646\unitlength}\raggedright $.4$\end{minipage}}}%
    \put(0.58161655,0.07910203){\color[rgb]{0,0,0}\makebox(0,0)[lt]{\begin{minipage}{0.0637168\unitlength}\raggedright $.6$\end{minipage}}}%
    \put(0.73717388,0.07910381){\color[rgb]{0,0,0}\makebox(0,0)[lt]{\begin{minipage}{0.08628319\unitlength}\raggedright $.8$\end{minipage}}}%
    \put(0.89272968,0.07910288){\color[rgb]{0,0,0}\makebox(0,0)[lt]{\begin{minipage}{0.07168138\unitlength}\raggedright $1$\end{minipage}}}%
    \put(0.10383898,0.14993663){\color[rgb]{0,0,0}\makebox(0,0)[lt]{\begin{minipage}{0.08893805\unitlength}\raggedright $0$\end{minipage}}}%
    \put(0.0982839,0.24993807){\color[rgb]{0,0,0}\makebox(0,0)[lt]{\begin{minipage}{0.08893804\unitlength}\raggedright $.2$\end{minipage}}}%
    \put(0.45513795,0.75238497){\color[rgb]{0,0,0}\makebox(0,0)[lt]{\begin{minipage}{0.37433627\unitlength}\raggedright $u(x,t)$\end{minipage}}}%
    \put(0.53284199,0.50585106){\color[rgb]{0,0,0}\makebox(0,0)[lt]{\begin{minipage}{0.16327429\unitlength}\raggedright $t=.1$\end{minipage}}}%
    \put(0.09828549,0.35271504){\color[rgb]{0,0,0}\makebox(0,0)[lt]{\begin{minipage}{0.08362832\unitlength}\raggedright $.4$\end{minipage}}}%
    \put(0.09828343,0.45271439){\color[rgb]{0,0,0}\makebox(0,0)[lt]{\begin{minipage}{0.08893805\unitlength}\raggedright $.6$\end{minipage}}}%
    \put(0.0982839,0.55549362){\color[rgb]{0,0,0}\makebox(0,0)[lt]{\begin{minipage}{0.12345132\unitlength}\raggedright $.8$\end{minipage}}}%
    \put(0.10661882,0.66382614){\color[rgb]{0,0,0}\makebox(0,0)[lt]{\begin{minipage}{0.08362832\unitlength}\raggedright $1$\end{minipage}}}%
    \put(0.43328447,0.3789828){\color[rgb]{0,0,0}\makebox(0,0)[lt]{\begin{minipage}{0.16725663\unitlength}\raggedright $t=.3$\end{minipage}}}%
  \end{picture}%
\endgroup%

%% file: ExF_paper.pdf_tex
\begingroup%
  \makeatletter%
  \providecommand\color[2][]{%
    \errmessage{(Inkscape) Color is used for the text in Inkscape, but the package 'color.sty' is not loaded}%
    \renewcommand\color[2][]{}%
  }%
  \providecommand\transparent[1]{%
    \errmessage{(Inkscape) Transparency is used (non-zero) for the text in Inkscape, but the package 'transparent.sty' is not loaded}%
    \renewcommand\transparent[1]{}%
  }%
  \providecommand\rotatebox[2]{#2}%
  \ifx\svgwidth\undefined%
    \setlength{\unitlength}{576.00001843bp}%
    \ifx\svgscale\undefined%
      \relax%
    \else%
      \setlength{\unitlength}{\unitlength * \real{\svgscale}}%
    \fi%
  \else%
    \setlength{\unitlength}{\svgwidth}%
  \fi%
  \global\let\svgwidth\undefined%
  \global\let\svgscale\undefined%
  \makeatother%
  \begin{picture}(1,0.75)%
    \put(0,0){\includegraphics[width=\unitlength,page=1]{ExF_paper.pdf}}%
    \put(0.86152929,0.68758654){\color[rgb]{0,0,0}\makebox(0,0)[lt]{\begin{minipage}{0.11814159\unitlength}\raggedright $t=.1$\end{minipage}}}%
    \put(0.65833363,0.36679009){\color[rgb]{0,0,0}\makebox(0,0)[lt]{\begin{minipage}{0.16725663\unitlength}\raggedright $t=5$\end{minipage}}}%
    \put(0.48180957,0.03279306){\color[rgb]{0,0,0}\makebox(0,0)[lt]{\begin{minipage}{0.28407079\unitlength}\raggedright $x$\end{minipage}}}%
    \put(0.12328366,0.07674092){\color[rgb]{0,0,0}\makebox(0,0)[lt]{\begin{minipage}{0.06106193\unitlength}\raggedright $0$\end{minipage}}}%
    \put(0.2705057,0.07677398){\color[rgb]{0,0,0}\makebox(0,0)[lt]{\begin{minipage}{0.07035399\unitlength}\raggedright $.2$\end{minipage}}}%
    \put(0.42606341,0.07910372){\color[rgb]{0,0,0}\makebox(0,0)[lt]{\begin{minipage}{0.079646\unitlength}\raggedright $.4$\end{minipage}}}%
    \put(0.58161655,0.07910203){\color[rgb]{0,0,0}\makebox(0,0)[lt]{\begin{minipage}{0.0637168\unitlength}\raggedright $.6$\end{minipage}}}%
    \put(0.73717388,0.07910381){\color[rgb]{0,0,0}\makebox(0,0)[lt]{\begin{minipage}{0.08628319\unitlength}\raggedright $.8$\end{minipage}}}%
    \put(0.89272968,0.07910288){\color[rgb]{0,0,0}\makebox(0,0)[lt]{\begin{minipage}{0.07168138\unitlength}\raggedright $1$\end{minipage}}}%
    \put(0.09828343,0.10271441){\color[rgb]{0,0,0}\makebox(0,0)[lt]{\begin{minipage}{0.08893805\unitlength}\raggedright $0$\end{minipage}}}%
    \put(0.0982839,0.23604918){\color[rgb]{0,0,0}\makebox(0,0)[lt]{\begin{minipage}{0.08893804\unitlength}\raggedright $.2$\end{minipage}}}%
    \put(0.45513795,0.75238497){\color[rgb]{0,0,0}\makebox(0,0)[lt]{\begin{minipage}{0.37433627\unitlength}\raggedright $u(x,t)$\end{minipage}}}%
    \put(0.74353519,0.60009884){\color[rgb]{0,0,0}\makebox(0,0)[lt]{\begin{minipage}{0.16327429\unitlength}\raggedright $t=.3$\end{minipage}}}%
    \put(0.09828549,0.37215948){\color[rgb]{0,0,0}\makebox(0,0)[lt]{\begin{minipage}{0.08362832\unitlength}\raggedright $.4$\end{minipage}}}%
    \put(0.09828343,0.50549217){\color[rgb]{0,0,0}\makebox(0,0)[lt]{\begin{minipage}{0.08893805\unitlength}\raggedright $.6$\end{minipage}}}%
    \put(0.0982839,0.64160473){\color[rgb]{0,0,0}\makebox(0,0)[lt]{\begin{minipage}{0.12345132\unitlength}\raggedright $.8$\end{minipage}}}%
  \end{picture}%
\endgroup%